\documentclass[11pt]{amsart}
\usepackage{}
\usepackage{amsmath}
\usepackage{amsfonts}
\usepackage{amssymb}
\usepackage[all,cmtip]{xy}           

\usepackage{dsfont}

\usepackage{bbm}
\usepackage{bbding}
\usepackage{txfonts}
\usepackage[shortlabels]{enumitem}
\usepackage{ifpdf}
\ifpdf
\usepackage[colorlinks,linkcolor=black,final,backref=page,hyperindex]{hyperref}
\else
\usepackage[colorlinks,linkcolor=black,final,backref=page,hyperindex,hypertex]{hyperref}
\fi
\usepackage{amscd}
\usepackage{tikz-cd}
\usepackage[active]{srcltx}

\makeatletter

\newtheorem{Thm}{Theorem}[section]
\newtheorem{Prop}[Thm]{Proposition}
\newtheorem{Lem}[Thm]{Lemma}
\newtheorem{Cor}[Thm]{Corollary}
\newtheorem{Def}[Thm]{Definition}
 \theoremstyle{remark}
 \newtheorem{Rem}[Thm]{Remark}
\newtheorem{exam}[Thm]{Example}

\newcommand{\nc}{\newcommand}

\newenvironment{Proof}[1][Proof]{\begin{trivlist}
\item[\hskip \labelsep {\bfseries #1}]}{\flushright
$\Box$\end{trivlist}}

\nc{\delete}[1]{{}}
\nc{\mmargin}[1]{}

\nc{\mlabel}[1]{\label{#1}}  
\nc{\mcite}[1]{\cite{#1}}  
\nc{\mref}[1]{\ref{#1}}  
\nc{\mbibitem}[1]{\bibitem{#1}} 

\delete{
	\nc{\mlabel}[1]{\label{#1}  
		{\hfill \hspace{1cm}{\bf{{\ }\hfill(#1)}}}}
	\nc{\mcite}[1]{\cite{#1}{{\bf{{\ }(#1)}}}}  
	\nc{\mref}[1]{\ref{#1}{{\bf{{\ }(#1)}}}}  
	\nc{\mbibitem}[1]{\bibitem[\bf #1]{#1}} 
}

 \font\cyrs=wncyr7
\newcommand{\Z}{{\mathbb{Z}}}

\newcommand{\M}{{\mathcal{M}}}
 \newcommand{\Path}{{\mathcal{P}}}
  \newcommand{\V}{{\mathcal{V}}}

  \newcommand{\CH}{{\mathrm{CH}}}

  \newcommand{\Tip}{{\mathrm{Tip}}}
 \newcommand{\Tor}{{\mathrm{Tor}}}
  \newcommand{\NonTip}{{\mathrm{NonTip}}}

\nc{\la}{\longrightarrow}
\nc{\ot}{\otimes}
\nc{\rar}{\rightarrow}
\nc{\lon }{\,\rightarrow\,}
\nc{\dar}{\downarrow}
\nc{\dap}[1]{\downarrow \rlap{$\scriptstyle{#1}$}}
\nc{\defeq}{\stackrel{\rm def}{=}}
\nc{\dis}[1]{\displaystyle{#1}}
\nc{\dotcup}{\ \displaystyle{\bigcup^\bullet}\ }
\nc{\hcm}{\ \hat{,}\ }
\nc{\hts}{\hat{\otimes}}
\nc{\hcirc}{\hat{\circ}}
\nc{\lleft}{[}
\nc{\lright}{]}
\nc{\curlyl}{\left \{ \begin{array}{c} {} \\ {} \end{array}
	\right .  \!\!\!\!\!\!\!}
\nc{\curlyr}{ \!\!\!\!\!\!\!
	\left . \begin{array}{c} {} \\ {} \end{array}
	\right \} }
\nc{\longmid}{\left | \begin{array}{c} {} \\ {} \end{array}
	\right . \!\!\!\!\!\!\!}
\nc{\ora}[1]{\stackrel{#1}{\rar}}
\nc{\ola}[1]{\stackrel{#1}{\la}}
\nc{\scs}[1]{\scriptstyle{#1}} \nc{\mrm}[1]{{\rm #1}}
\nc{\dirlim}{\displaystyle{\lim_{\longrightarrow}}\,}
\nc{\invlim}{\displaystyle{\lim_{\longleftarrow}}\,}
\nc{\dislim}[1]{\displaystyle{\lim_{#1}}} \nc{\colim}{\mrm{colim}}
\nc{\mvp}{\vspace{0.3cm}} \nc{\tk}{^{(k)}} \nc{\tp}{^\prime}
\nc{\ttp}{^{\prime\prime}} \nc{\svp}{\vspace{2cm}}
\nc{\vp}{\vspace{8cm}}
\nc{\modg}[1]{\!<\!\!{#1}\!\!>}
\nc{\intg}[1]{F_C(#1)}
\nc{\lmodg}{\!<\!\!}
\nc{\rmodg}{\!\!>\!}
\nc{\cpi}{\widehat{\Pi}}
\nc{\ssha}{{\mbox{\cyrs X}}} 
\nc{\tsha}{{\mbox{\cyrt X}}}
\nc{\shpr}{\diamond}    
\nc{\labs}{\mid\!}
\nc{\rabs}{\!\mid}

\nc{\gldim}{\mrm{gldim}}
\nc{\ad}{\mrm{ad}}
\nc{\ann}{\mrm{ann}}
\nc{\Aut}{\mrm{Aut}}
\nc{\bim}{\mbox{-}\mathsf{Bimod}}
\nc{\br}{\mrm{bre}}
\nc{\can}{\mrm{can}}
\nc{\Cont}{\mrm{Cont}}
\nc{\rchar}{\mrm{char}}
\nc{\cok}{\mrm{coker}}
\nc{\de}{\mrm{dep}}
\nc{\dtf}{{R-{\rm tf}}}
\nc{\dtor}{{R-{\rm tor}}}

\nc{\Div}{{\mrm Div}}
\nc{\Diff}{\mrm{DA}}
\nc{\Diffl}{\mathsf{DA}_\lambda}
\nc{\diffo}{{\mathsf{DO}_\lambda}}
\nc{\alg}{\mathsf{Alg}}
\nc{\End}{\mrm{End}}
\nc{\Ext}{\mrm{Ext}}
\nc{\Fil}{\mrm{Fil}}
\nc{\Fr}{\mrm{Fr}}
\nc{\Frob}{\mrm{Frob}}
\nc{\Gal}{\mrm{Gal}}
\nc{\GL}{\mrm{GL}}
\nc{\Hom}{\mrm{Hom}}
\nc{\Hoch}{\mrm{Hoch}}
\nc{\hsr}{\mrm{H}}
\nc{\hpol}{\mrm{HP}}
\nc{\id}{\mrm{id}}
\nc{\im}{\mrm{im}}
\nc{\Id}{\mrm{Id}}
\nc{\ID}{\mrm{ID}}
\nc{\Irr}{\mrm{Irr}}
\nc{\incl}{\mrm{incl}}
\nc{\length}{\mrm{length}}
\nc{\NLSW}{\mrm{NLSW}}
\nc{\Lie}{\mrm{Lie}}
\nc{\mchar}{\rm char}
\nc{\mpart}{\mrm{part}}
\nc{\ql}{{\QQ_\ell}}
\nc{\qp}{{\QQ_p}}
\nc{\rank}{\mrm{rank}}
\nc{\rcot}{\mrm{cot}}
\nc{\rdef}{\mrm{def}}
\nc{\rdiv}{{\rm div}}
\nc{\rtf}{{\rm tf}}
\nc{\rtor}{{\rm tor}}
\nc{\res}{\mrm{res}}
\nc{\SL}{\mrm{SL}}
\nc{\Spec}{\mrm{Spec}}
\nc{\tor}{\mrm{tor}}
\nc{\Tr}{\mrm{Tr}}
\nc{\tr}{\mrm{tr}}
\nc{\wt}{\mrm{wt}}

\nc{\bfk}{{\bf k}}
\nc{\bfone}{{\bf 1}}
\nc{\bfzero}{{\bf 0}}
\nc{\detail}{\marginpar{\bf More detail}
	\noindent{\bf Need more detail!}
	\svp}
\nc{\gap}{\marginpar{\bf Incomplete}\noindent{\bf Incomplete!!}
	\svp}
\nc{\FMod}{\mathbf{FMod}}
\nc{\Int}{\mathbf{Int}}
\nc{\Mon}{\mathbf{Mon}}
\nc{\remarks}{\noindent{\bf Remarks: }}
\nc{\Rep}{\mathbf{Rep}}
\nc{\Rings}{\mathbf{Rings}}
\nc{\Sets}{\mathbf{Sets}}
\nc{\ob}{\mathsf{Ob}}

\nc{\BA}{{\mathbb A}}   \nc{\CC}{{\mathbb C}}
\nc{\DD}{{\mathbb D}}   \nc{\EE}{{\mathbb E}}
\nc{\FF}{{\mathbb F}}   \nc{\GG}{{\mathbb G}}
\nc{\HH}{{\mathbb H}}   \nc{\LL}{{\mathbb L}}
\nc{\NN}{{\mathbb N}}   \nc{\PP}{{\mathbb P}}
\nc{\QQ}{{\mathbb Q}}   \nc{\RR}{{\mathbb R}}
\nc{\TT}{{\mathbb T}}   \nc{\VV}{{\mathbb V}}
\nc{\ZZ}{{\mathbb Z}}   \nc{\TP}{\widetilde{P}}
\nc{\m}{{\mathbbm m}}

\nc{\cala}{{\mathcal A}}    \nc{\calc}{{\mathcal C}}
\nc{\cald}{\mathcal{D}}     \nc{\cale}{{\mathcal E}}
\nc{\calf}{{\mathcal F}}    \nc{\calg}{{\mathcal G}}
\nc{\calh}{{\mathcal H}}    \nc{\cali}{{\mathcal I}}
\nc{\call}{{\mathcal L}}    \nc{\calm}{{\mathcal M}}
\nc{\caln}{{\mathcal N}}    \nc{\calo}{{\mathcal O}}
\nc{\calp}{{\mathcal P}}    \nc{\calr}{{\mathcal R}}
\nc{\cals}{{\mathcal S}}    \nc{\calt}{{\Omega}}
\nc{\calv}{{\mathcal V}}    \nc{\calw}{{\mathcal W}}
\nc{\calx}{{\mathcal X}}
\nc{\U}{{\mathcal U}} \nc{\D}{{\mathcal D}}

\nc{\fraka}{{\mathfrak a}}
\nc{\frakb}{\mathfrak{b}}
\nc{\frakg}{{\frak g}}
\nc{\frakl}{{\frak l}}
\nc{\fraks}{{\frak s}}
\nc{\frakB}{{\frak B}}
\nc{\frakm}{{\frak m}}
\nc{\frakM}{{\frak M}}
\nc{\frakp}{{\frak p}}
\nc{\frakW}{{\frak W}}
\nc{\frakX}{{\frak X}}
\nc{\frakS}{{\frak S}}
\nc{\frakA}{{\frak A}}
\nc{\frakx}{{\frakx}}

\begin{document}

\title[algebraic Morse theory]{
Algebraic Morse theory via Homological Perturbation Lemma}

\author{Jun Chen}

\address{Jun Chen, School of Mathematics,
   Nanjing University,
   Nanjing 210093,
   P.R.China}
\email{mathcj@nju.edu.cn}

\author{Yuming Liu}
\address{Yuming Liu,
  School of Mathematical Sciences,
  Laboratory of Mathematics and Complex Systems,
  Beijing Normal University,
  Beijing 100875,
  P.R.China}
\email{ymliu@bnu.edu.cn}

\author{Guodong Zhou}

\address{Guodong Zhou,  School of Mathematical Sciences, Key Laboratory of MEA (Ministry of Education), Shanghai Key Laboratory of PMMP,
  East China Normal University,
 Shanghai 200241,
   P.R.China}

\email{gdzhou@math.ecnu.edu.cn}

\date{\today}

\begin{abstract}
   As a  generalization of  the classical killing-contractible-complexes lemma,  we present   algebraic Morse theory  via homological perturbation lemma,  in a   form more general than existing presentations in the literature.   Two-sided Anick resolutions due to  E.~Sk\"{o}ldberg are generalised to algebras given by quivers with relations and    a  minimality criterion is provided as well.     Two applications of algebraic Morse theory are presented.  It is shown that the Chinese algebra of rank $n\geq 1$ is homologically smooth and of global dimension  $\frac{n(n+1)}{2}$, and the minimal two-sided projective resolution of a Koszul algebra is constructed.
\end{abstract}

\subjclass[2010]{18G35, 16E05, 16E10}

\keywords{Algebraic Morse theory; Chinese algebra; Homological perturbation lemma; Killing-contractible-complexes lemma; Two-sided Anick resolution}

\maketitle

 \tableofcontents

\allowdisplaybreaks

\section*{Introduction}

Discrete Morse theory has  its origin in the work of K.~S.~Brown, R.~Geoghegan \cite{BG}. In that paper, one encountered  a cell complex   with one
vertex and infinitely many cells in each positive dimension. The authors, using ad hoc method,
 collapsed this cell complex to a quotient complex with only two cells in each
positive dimension. K.~S.~Brown formalized  and
applied the  collapsing method  scheme to groups with  a rewriting system \cite{Brown}. Motivated by differential topology, R.~Forman \cite{Forman} rediscovered and developed this  theory  as a discrete version of the usual smooth Morse  theory, hence the name ``discrete Morse theory".
 Since then, this subject has received much attention in combinatorial and computational topology; see  \cite{Forman02, Kozlov08, Knudson, Scoville} etc.

An algebraic version of discrete Morse theory has been developed by E.~Sk\"{o}ldberg \cite{S06}, D.~Kozlov \cite{Kozlov}, M.~J\"{o}llenbeck and V.~Welker \cite{JW}.  It has many applications in algebra, such as
combinatorial commutative algebra \cite{JW}, cohomology of Lie algebras \cite{S05, LV16, LV17b, LV19}, Hochschild cohomology \cite{Lopatkin, LV17a}, and operad  theory \cite{DK13}.

E.~Sk\"{o}ldberg \cite{S18} studied  algebraic Morse theory from the viewpoint of the Homological Perturbation Lemma \cite{EM, Shih, Brown64, HK, Crainic}.
This paper pushes this idea further, generalising and complementing   \cite{S06, S18} as well as \cite{Kozlov, JW} and presents applications to Chinese algebras and a Koszul algebra.

\medskip

The layout of this paper is as follows.

The first section contains an introduction to the Homological Perturbation Lemma (HPL) by providing some new proofs and   observations which seem to be of independent interest.
In the second section we use   HPL  to generalise   the classical killing-contractible-complexes lemma in homological algebra.

Algebraic Morse theory is presented in the third section  and our definition of the key concept ``Morse matching" is  more general than those in
\cite{S06,Kozlov,JW}.

One of the most important advantages of algebraic Morse theory is that it provides a method to construct two-sided Anick resolutions, compared with the original one-sided version in \cite{Anick, AG}. This is particularly useful for Hochschild cohomology.
In the fourth section, using the corresponding Gr\"{o}bner-Shirshov basis theory \cite{Green}, we  furnish a construction of two-sided Anick resolutions   for  algebras given by quivers with relations,   generalising  \cite[Section 3.2]{S06}.
A similar resolution with a different construction of the differentials is provided by  S.~Chouhy and A.~Solotar \cite{CS}, which might be closely related to two-sided Anick resolutions.

 The fifth section contains a   criterion  which gives a sufficient condition when  the two-sided Anick resolution  is minimal;  see Theorem~\ref{Thm: minimal criterion}. We also provide a counter-example to    a result that appeared in \cite{JW}.

The sixth  section contains an application of algebraic Morse theory to Chinese algebras \cite{DK92,CEKNH}. Using the
Gr\"{o}bner-Shirshov basis  for  Chinese algebras discovered by Y.~Chen and J.~Qiu \cite{CQ}, we prove that  a Chinese algebra  of   rank $n\geq 1$ is  homologically smooth, that is, it admits a finite length resolution by finitely generated projective bimodules,  and  its  global  dimension is     $\frac{n(n+1)}{2}$; see Theorem~\ref{Thm: Chinese algebras}.

 In last section, we consider the quadratic algebra $A=k \langle x,y,z \ |\  x^2+yx ,xz,zy\rangle$ first studied by  N.~Iyudu and S.~Shkarin  \cite{IS} in their classification result for Hilbert series of Koszul algebras with three generators and three relations with the goal to answer two questions from the book \cite{PP}. In \cite{Dots}, V.~Dotsenko and S.~R.~Chowdhury constructed the one-sided Anick resolution of $A$. We will use the algebraic Morse theory to construct the two-sided minimal resolution of $A$ through its two-sided Anick resolution.

\bigskip

\section{Homological Perturbation Lemma}

 In this section, we  give  an  introduction to the Homological Perturbation Lemma (HPL).

We first introduce some relevant notions  which appear to be well known, although to our knowledge,   some of them  do not  appear  in the literature.
 \begin{Def}
 \begin{itemize}
 \item[(a)]  A  homotopy  retract datum  (aka. HR datum)
 \begin{equation}\xymatrix{{(L_*,b)}\ar@<1ex>[r]^{i} & {(M_*,b)} \ar@<1ex>[l]^{p}&\ar@(ur,dr)^{h}}    \label{equ: HR data} \end{equation}
consists of the following
\begin{itemize}
\item[(i)] two chain maps $i$ and $p$ between  complexes $(L_*, b)$ and $(M_*, b)$,  and
\item[(ii)] a homotopy $h$ between $ip$ and $1$ (so $ip=1+bh+hb$).
\end{itemize}
 \item[(b)] An  HR datum   is a strong quasi-isomorphism datum (aka. SQI datum) if moreover,
\begin{itemize}
\item[(iii)]  $p$ and $i$ are  quasi-isomorphisms.
\end{itemize}
\item[(c)] An  SQI datum   is a homotopy equivalence datum  (aka. HE datum) if moreover,
\begin{itemize}
\item[(iv)]  $pi$ is homotopic to $1$.
\end{itemize}
\item[(d)]  An HE datum   is a deformation retract datum  (aka. DR datum),  if moreover,
\begin{itemize}
\item[(v)]  $pi=1$.
\end{itemize}
\item[(e)] A  DR datum is a   strong  deformation retract datum (aka. SDR datum),  if moreover,
\begin{itemize}
\item[(vi)]  $h^2=0, hi=0, ph=0$.
\end{itemize}
\end{itemize}
\end{Def}

 \begin{Rem}\begin{itemize}
  \item[(a)]   DR data and SDR data exist even from  the beginning of HPL \cite{EM,    Shih, Brown64}.
  The name HE data appeared in the work of M.~Crainic \cite{Crainic} with a different meaning, namely that which we call SQI data.

\item[(b)]   A DR datum can be always modified to an SDR datum, as explained in    \cite{LS}.

\end{itemize}

\end{Rem}

We now introduce  a condition  which lies at the heart of the HPL.

 A perturbation $\delta$ of Datum~(\ref{equ: HR data}) is a graded  map on $M_*$ of the same degree as $b$ such that $b+\delta$ is a new differential. We call it \textit{small}, if $1-\delta h$ is invertible. In this case, put
$$A=(1-\delta h)^{-1} \delta$$
and consider a new perturbed datum
  \begin{equation}\xymatrix{(L_*, b_\infty)\ar@<1ex>[r]^{i_\infty} & (M_*, b+\delta )    \ar@<1ex>[l]^{p_\infty}  & \ar@(ur, dr)^{h_\infty}  }   \label{equ: new perturbed data} \end{equation}
with
$$i_\infty=i+hAi, p_\infty=p+pAh, h_\infty=h+hAh, b_\infty=b+pAi.  $$

\begin{Rem} In a heuristic manner,
$$\begin{array}{rclcl} A&=&(1-\delta h)^{-1} \delta&=&\sum_{n=0}^{\infty} (\delta h)^n\delta=\sum_{n=0}^{\infty} \delta(h\delta)^n,\\
 i_\infty&=&i+hAi&=&\sum_{n=0}^\infty (h\delta)^n i,\\
  p_\infty&=&p+pAh&=&\sum_{n=0}^{\infty} p(\delta h)^n,\\
   h_\infty&=&h+hAh&=&\sum_{n=0}^{\infty} h(\delta h)^n,\\
   b_\infty&=&b+pAi&=& b+\sum_{n=0}^{\infty} p(\delta h)^n\delta i=b+\sum_{n=0}^{\infty} p\delta(h\delta)^n i.\end{array}  $$

\end{Rem}

The following result summarizes known versions of the HPL.
\begin{Thm}[Homological Perturbation Lemma]\label{Thm: perturbation lemma}
 Assume that  $\delta$ is a small perturbation of Datum~(\ref{equ: HR data}). If  Datum~(\ref{equ: HR data}) is  an HR datum (resp. SQI datum, HE datum,   SDR datum), then so is
Datum~(\ref{equ: new perturbed data}).
\end{Thm}

\begin{Rem}

\begin{itemize}
\item[(a)]
  It is well known that Theorem~\ref{Thm: perturbation lemma} does NOT hold for DR data; see for example,   \cite[2.3 Remarks (i)]{Crainic}.

  \item[(b)] The SDR version is the first version of the HPL that appeared in the literature \cite{EM, Shih, Brown64};   the HR version is essentially contained in the proof of the SDR version, in particular, we have  $$i_\infty p_{\infty}=1+(b+\delta)h_{\infty}+h_{\infty}(b+\delta);$$  the HE version was first shown in \cite{HK}; the SQI version is contained in \cite{Crainic}.   For a   detailed historic account, see \cite{Hue}.

\end{itemize}
\end{Rem}

We will present a direct  proof for the HE version compared with the original proof of \cite{HK}, which used a mapping cylinder construction to deduce  the HE version from the SDR version,  and we will also provide a streamlined proof for the SQI version. In the course of the proofs, we add some which appear to be new observations.

\medskip

In the proofs of this section,  in order to facilitate the understanding of the reader, we often underline some terms indicating that the terms  will be changed in the next equality.

For the rest of this section, assume that we are given an HR datum (\ref{equ: HR data}) which is
 endowed with a small perturbation $\delta$.

The following three equalities are crucial in the proof  of  the  HR version of HPL.
\begin{Lem}[{see \cite[Main Perturbation Lemma]{Crainic}}] \label{Lem: several equations}
 We have the following equalities:
 \begin{equation} \label{Eq: Ah8}
 \delta hA=Ah\delta=A-\delta ,
 \end{equation}
 \begin{equation} \label{Eq: 1+Ah}
 (1- \delta h )^{-1}= 1+ Ah , (1-h\delta)^{-1} =1 +h A ,
 \end{equation}
 \begin{equation}\label{Eq: AipA}
 AipA+Ab+bA=0.
 \end{equation}
\end{Lem}



 From now on we assume both  the  HR version and  the  SDR version of HPL. We are going to deduce the SQI and  the  HE version of HPL from these assumptions.
The following  interesting observations  are useful in the proof of the  SQI and the  HE version of  HPL.

\begin{Lem}\label{Lem: p_inf i-1 p i_inf-1}
 The maps
 $p_{\infty} i -1 : (L_*, b) \longrightarrow (L_*, b_{\infty})$
and
 $p i_{\infty} -1 : (L_*, b_{\infty}) \longrightarrow (L_*, b)$
are chain maps, that is,
 \begin{equation} \label{Eq: pinf i-1}
b_{\infty} (p_{\infty}i-1)=(p_{\infty}i-1)b
\end{equation}
\begin{equation} \label{Eq: p iinf-1}
(p i_{\infty}-1)b_{\infty}=b(p i_{\infty}-1).
\end{equation}

\end{Lem}

\begin{Proof} We will show   Equality~(\ref{Eq: pinf i-1}) and Equality~(\ref{Eq: p iinf-1}) can be proved similarly.
In fact, we have
$$\begin{array}{rcl}
b_{\infty} (p_{\infty}i-1) &=&\underline{ b_{\infty} p_{\infty}}i - b_{\infty}\\
&= &p_{\infty}  bi+\underline{p_{\infty}} \delta i   -b_{\infty} \\
&= &p_{\infty}  bi+p\underline{(\delta+Ah\delta)}   i   -b_{\infty} \\
 &\stackrel{(\ref{Eq: Ah8})}{=} &p_{\infty}\underline{ b i}+ \underline{pAi-b_{\infty}}\\
 &=& p_{\infty} i b- b\\
 &=& (p_{\infty} i -1)b.
\end{array}$$
\end{Proof}



\begin{Lem}\label{Lem: nullhomotoyp i(pi_infty-1)}
 The chain maps
$$i \circ (p i_{\infty} - 1 ) : (L_*, b_{\infty})\longrightarrow (M_*, b),$$
$$(p_{\infty} i - 1 ) \circ p : (M_*, b) \longrightarrow (L_*, b_{\infty}),$$
$$(p i_{\infty} - 1 ) \circ p_{\infty} : (M_*, b+\delta) \longrightarrow (L_*, b),$$
and
$$i_{\infty} \circ (p_{\infty} i - 1 ) : (L_*, b) \longrightarrow (M_*, b+\delta)$$
are all null-homotopic.

\end{Lem}

\begin{Proof}
In fact,     we have the following equalities:
\begin{equation}\label{Eq: i(pi_inf-1)} i(pi_{\infty}-1) = hi_{\infty} b_{\infty}+ bhi_{\infty},\end{equation}
\begin{equation}\label{Eq: (p_inf i-1)p} (p_{\infty} i -1)p = b_{\infty} p_{\infty} h + p_{\infty} h b, \end{equation}
\begin{equation}\label{Eq: (pi_inf -1)p_inf}  (pi_{\infty}-1)p_{\infty}= p h_{\infty}(b+\delta) + b p h_{\infty},  \end{equation}
\begin{equation}\label{Eq: i_inf(p_inf i-1)} i_{\infty}(p_{\infty} i -1) =(b+\delta )h_{\infty} i + h_{\infty} i b. \end{equation}
  whose proofs can be  verified directly.

%

\end{Proof}

\begin{Lem}\label{Lem: pi-1 and p_inf i_inf -1}
Let  $h'=ph_{\infty}i$ and $h''=p_{\infty} h i_{\infty}$. The following equalities hold:
\begin{equation} \label{Eq: pi-1} pi-1=b h' + h'b-(pi_{\infty}-1)(p_{\infty} i-1),\end{equation}
 and
 \begin{equation} \label{Eq:  p_infty i_infty -1} p_{\infty} i_{\infty} -1 = b_{\infty} h'' + h'' b_{\infty} -(p_{\infty}i -1)(p i_{\infty} - 1).\end{equation}

\end{Lem}

\begin{Proof}  We only prove Equality~(\ref{Eq:  p_infty i_infty -1}), as Equality~(\ref{Eq: pi-1}) can be shown similarly.


%

 For Equality~(\ref{Eq:  p_infty i_infty -1}),
   $$\begin{array}{rcl}
   p_{\infty} i_{\infty} -1 - b_{\infty} h'' - h'' b_{\infty} &=& p_{\infty} i_{\infty} -1 - \underline{b_{\infty} p_{\infty}} h i_{\infty} - p_{\infty} h \underline{i_{\infty}b_{\infty}}\\
   & = &  p_{\infty} i_{\infty} -1 -  p_{\infty} ( \underline{b} + \delta )\underline{h} i_{\infty} - p_{\infty} \underline{h }(\underline{b}+\delta) i_{\infty}\\
   &=& \underline{p_{\infty} i_{\infty}} -1 - p_{\infty} ip i_{\infty} + \underline{p_{\infty} i_{\infty}} -\underline{ p_{\infty}} \delta h i_{\infty} - p_{\infty} h \delta \underline{i_{\infty}}\\
   &=& - p_{\infty} ip i_{\infty} + 2 p_{\infty} i_{\infty} -1 - p \underline{(1+Ah) \delta h} i_{\infty} -p_{\infty} h \underline{\delta (1+hA)} i\\
   &\overset{(\ref{Eq: Ah8})}{=}& - p_{\infty} ip i_{\infty} + \underline{2 p_{\infty} i_{\infty}} -1 - \underline{p Ah i_{\infty} } -\underline{p_{\infty} h A i}\\
   &=& - p_{\infty} ip i_{\infty} - 1 +  p i_{\infty} +p_{\infty} i\\
   &=& -(p_{\infty}i-1)(p i_{\infty} - 1),
   \end{array}$$
   where the second equality follows from  the fact that $p_\infty$ and $i_\infty$ are chain maps (by the HR version) .
\end{Proof}

We now prove the SQI and  the  HE version of Theorem~\ref{Thm: perturbation lemma}.

\smallskip

\begin{Proof}[Proof of  the   SQI  version of  Theorem~\ref{Thm: perturbation lemma} ]\
 Since an SQI datum~(\ref{equ: HR data})   inherently satisfies  the conditions of an HR datum,   the perturbed datum~(\ref{equ: new perturbed data}) is at least an HR datum.
This implies that the composition  $i_{\infty} p_{\infty} $ induces the identity map on homology.
Furthermore,  the chain map
 $$ b_{\infty} h'' + h'' b_{\infty}:(L_*, b_{\infty}) \to (L_*, b_{\infty})   $$
 vanishes  on  homology. By equation~\eqref{Eq:  p_infty i_infty -1}, the proof now reduces to demonstrating that $(p_{\infty}i -1)(p i_{\infty} - 1)$ vanishes on homology.

 Lemma~\ref{Lem: nullhomotoyp i(pi_infty-1)} establishes that the chain map $i  (p i_{\infty} - 1 )$ is null-homotopic, hence induces the zero map on homology. Since  $i$ is a quasi-isomorphism, it follows that
 $$ (p i_{\infty} - 1 ):  (L_*, b_{\infty}) \to (L_*, b) $$
must induce the zero map on homology. By the functoriality of homology, the composition $(p_{\infty}i -1)(p i_{\infty} - 1) $ consequently vanishes on homology.
 \end{Proof}

\begin{Proof}[Proof of the HE version of  Theorem~\ref{Thm: perturbation lemma}]\

Given a small perturbation $\delta$ on an HE datum~(\ref{equ: HR data}) with a homotopy $k$ on $(L_*, b)$ such that $$pi-1=bk+kb$$
we need to construct a new homotopy $k_{\infty}$ on $(L_*, b_{\infty})$ such that $$p_{\infty} i_{\infty} - 1 =b_{\infty} k_{\infty} + k_{\infty} b_{\infty}.$$

Write $k'= (p_{\infty}i -1) k (p i_{\infty}-1) $.
We obtain
$$\begin{array}{rcl}
&&b_{\infty} k' + k' b_{\infty}\\
 & = & \underline{b_{\infty}(p_{\infty}i -1)} k (p i_{\infty}-1) + (p_{\infty}i -1) k \underline{(p i_{\infty}-1)  b_{\infty}}\\
&\overset{(\ref{Eq: pinf i-1})(\ref{Eq: p iinf-1})}{=}& (p_{\infty}i -1) \underline{bk} (p i_{\infty}-1) + (p_{\infty}i -1)\underline{ k b}(p i_{\infty}-1)\\
&=& (p_{\infty}i -1) (pi-1) (p i_{\infty}-1)\\
&=& \underline{(p_{\infty}i -1) p}~\underline{i (p i_{\infty}-1)} - (p_{\infty}i -1) (p i_{\infty}-1)\\
&\overset{(\ref{Eq: i(pi_inf-1)})(\ref{Eq: (p_inf i-1)p})}{=}& (b_{\infty} p_{\infty}h + p_{\infty} h b) (hi_{\infty} b_{\infty} + b h i_{\infty})- (p_{\infty}i -1) (p i_{\infty}-1)\\
&=& b_{\infty} p_{\infty} hh i_{\infty} b_{\infty} +p_{\infty} hh i_{\infty} b_{\infty} b_{\infty}  + b_{\infty} p_{\infty} h b h i _{\infty} + p_{\infty} hbhi_{\infty} b_{\infty}  - (p_{\infty}i -1) (p i_{\infty}-1)\\
&\overset{(\ref{Eq:  p_infty i_infty -1})}{=}& b_{\infty} p_{\infty} hh i_{\infty} b_{\infty} +p_{\infty} hh i_{\infty} b_{\infty} b_{\infty}
  + b_{\infty} p_{\infty} h b h i _{\infty} + p_{\infty} hbhi_{\infty} b_{\infty}- b_{\infty} p_{\infty} h i_{\infty} -p_{\infty} h i_{\infty} b_{\infty}
 \\
 &&+p_{\infty} i_{\infty} -1 .
\end{array}$$
Denote  $$k_{\infty} = (p_{\infty}i -1) k (p i_{\infty}-1)  +p_{\infty} h i_{\infty} -p_{\infty} hbh i_{\infty} -p_{\infty} hh i_{\infty} b_{\infty},$$
then we get the desired equality:
  $$p_{\infty}i_{\infty}-1=b_{\infty} k_{\infty} + k_{\infty} b_{\infty}.$$
\end{Proof}

\bigskip

\section{A generalization of killing-contractible-complexes lemma}

In  this section, given a DR datum~(\ref{equ: HR data}), we   describe the cokernel of $i_\infty$ in the perturbed datum (\ref{equ: new perturbed data}).  We will see that this is a generalization of the classical killing-contractible-complexes lemma \cite[Lemma 2.1.6]{Loday}.

\begin{Thm} \label{Thm: reformulation-algebraic-morse} Let $R$ be an associative ring, and let $$ \cdots\longrightarrow
 C_{n}\oplus C_{n}' \stackrel{d=\left(
\begin{array}{cc} \alpha&\beta \\ \gamma& \eta  \end{array}\right)}{\longrightarrow}
 C_{n-1}\oplus C_{n-1}' \longrightarrow\cdots$$ be a chain complex of $R$-modules. Assume that $\eta=\eta'+\eta''$ such that $(C_*',\eta')$ is a complex and is contractible with contracting homotopy $ \sigma  : C_n'\longrightarrow C_{n+1}'$. Suppose moreover, that $1+\eta''\sigma : C_n'\longrightarrow C_n'$ is  invertible (for example, $\eta''\sigma $ is  locally nilpotent, that is, for any $x\in C_n'$, there exists a positive integer $p$ such that $(\eta''\sigma )^p(x)=0$.). Denote the map $\lambda=(1+\eta''\sigma )^{-1}$.  Then the following statements hold:

\begin{itemize}\item [(a)]  $(C_*, \overline{d}:=\alpha-\beta \sigma \lambda\gamma)$  is a complex;




 \item[(b)]  there  exist  two short exact sequences of complexes
 $$0\to (C_*',\eta') \stackrel{\left(
\begin{array}{c} \beta \sigma  \\ 1+\eta''\sigma  \end{array}\right)}{\longrightarrow} (C_*\oplus C_*',d) \stackrel{g=(1\ \ -\beta \sigma \lambda)}{\longrightarrow} (C_*, \alpha-\beta \sigma \lambda\gamma)\to 0 $$and
 $$0\to (C_*,\alpha-\beta \sigma \lambda\gamma) \stackrel{f=\left(
\begin{array}{c} 1 \\ -\sigma \lambda\gamma \end{array}\right)}{\longrightarrow} (C_*\oplus C_*',d) \stackrel{(\sigma  \gamma\ \  1+\sigma \eta'')}{\longrightarrow} (C_*', \eta')\to 0;$$

\item[(c)]  $f$ and $g$  establish a homotopy equivalence between $(C_*, \overline{d})$ and $(C_{*}\oplus C_{*}', d)$.

  \end{itemize}

\end{Thm}

 \begin{Proof}  Let $$(L_*, b)=(C_*, 0), (M_*, b)=\Big(C_*\oplus C_*',\left(
\begin{array}{cc} 0 & 0 \\  0 &\eta' \end{array}\right) \Big),$$
and $$i=\left(
\begin{array}{c} 1 \\ 0 \end{array}\right), p=(1, 0), h=\left(
\begin{array}{cc} 0 &0 \\ 0 & -\sigma  \end{array}\right),
\delta=\left(
\begin{array}{cc} \alpha&\beta \\ \gamma&\eta'' \end{array}\right).$$ In this case, we have  $pi=1$, $ph=0$ and $hi=0$, but not necessarily $h^2=0$. So the original datum is a DR datum, but not an SDR datum.
Since
$$1-\delta h=\left(
\begin{array}{cc} 1 &0  \\ 0&1 \end{array}\right)- \left(
\begin{array}{cc} \alpha&\beta \\ \gamma&\eta'' \end{array}\right)\left(
\begin{array}{cc} 0 &0 \\ 0 & -\sigma  \end{array}\right)=\left(
\begin{array}{cc} 1 &\beta\sigma \\ 0 & 1+\eta''\sigma  \end{array}\right),$$
the invertibility of $1+\eta''\sigma$ implies that $\delta$ is small, and as a DR datum is always an HE datum,  we can apply the HE version of Theorem~\ref{Thm: perturbation lemma}.

\medskip

Let us compute the new perturbed datum.
We have $$A=(1-\delta h)^{-1} \delta =\left(
\begin{array}{cc} 1 &\beta\sigma \\ 0 & 1+\eta''\sigma  \end{array}\right)^{-1}\left(
\begin{array}{cc} \alpha&\beta \\ \gamma&\eta'' \end{array}\right)=\left(
 \begin{array}{cc} \alpha-\beta\sigma\lambda\gamma &\beta-\beta\sigma\lambda \eta'' \\ \lambda \gamma & \lambda \eta''  \end{array}\right);$$
the new differential on $C_*$ is
   $$\overline{d}:=b_\infty=b+pAi=0+(1\  0)\left(
\begin{array}{cc} \alpha-\beta\sigma\lambda\gamma &\beta-\beta\sigma\lambda \eta'' \\ \lambda \gamma & \lambda \eta''  \end{array}\right)\left( \begin{array}{c} 1 \\ 0 \end{array} \right)=\alpha-\beta \sigma \lambda \gamma;$$
   we also get two chain maps
   $$f:=i_\infty=i+hAi=\left( \begin{array}{c} 1 \\ 0 \end{array} \right)+\left( \begin{array}{cc} 0&0 \\ 0&{-\sigma} \end{array} \right)\left( \begin{array}{cc} \alpha-\beta \sigma\lambda \gamma &\ \beta-\beta \sigma\lambda \eta'' \\ \lambda \gamma & \lambda \eta''  \end{array} \right)\left( \begin{array}{c} 1 \\ 0 \end{array} \right)=\left( \begin{array}{cc} 1 \\ -\sigma \lambda \gamma \end{array} \right)$$ and
   $$g:=p_\infty=p+pAh= (1\  0)+ (1\  0)\left( \begin{array}{cc} \alpha-\beta \sigma\lambda \gamma &\ \beta-\beta \sigma\lambda \eta'' \\ \lambda \gamma & \lambda \eta''  \end{array} \right)\left( \begin{array}{cc} 0&0 \\ 0&{-\sigma} \end{array} \right) =
  \left( \begin{array}{cc} 1 &-\beta \sigma \lambda \end{array} \right);$$
  the map
  $$h_\infty=h+hAh=\left( \begin{array}{cc} 0&0 \\ 0&{-\sigma} \end{array} \right)+\left( \begin{array}{cc} 0&0 \\ 0&{-\sigma} \end{array} \right)\left( \begin{array}{cc} \alpha-\beta \sigma\lambda \gamma &\ \beta-\beta \sigma\lambda \eta'' \\ \lambda \gamma & \lambda \eta''  \end{array} \right)\left( \begin{array}{cc} 0&0 \\ 0&{-\sigma} \end{array} \right)= \left( \begin{array}{cc} 0 &0 \\ 0 & -\sigma\lambda \end{array} \right)$$ realised a homotopy between $fg$ and $1$, and the map
    $$\begin{array}{rcl}  k_\infty&=& (p_\infty i -1) k (pi_\infty -1)+p_\infty h i_\infty-p_\infty h b h i_\infty -p_\infty h^2i_\infty b_\infty\\
 &=& (p_\infty i -1) k (pi_\infty -1) + p_\infty(  h - h b h  -h^2 (b+\delta )) i_\infty\\
  &=& \left( \begin{array}{cc} 1 &-\beta \sigma \lambda \end{array} \right)
 \left[ \left(\begin{array}{cc} 0 &0 \\ 0 & -\sigma  \end{array}\right)
  - \left(\begin{array}{cc} 0 &0 \\ 0 & -\sigma  \end{array}\right)
 \left(\begin{array}{cc} 0 &0 \\ 0 & \eta'  \end{array}\right)
 \left(\begin{array}{cc} 0 &0 \\ 0 & -\sigma  \end{array}\right)
 -\left(\begin{array}{cc} 0 &0 \\ 0 & -\sigma  \end{array}\right)^2 \left(\begin{array}{cc} \alpha & \beta  \\ \gamma & \eta \end{array}\right) \right]
\left( \begin{array}{cc} 1 \\ -\sigma \lambda \gamma \end{array} \right)\\
&=& \left( \begin{array}{cc} 1 &-\beta \sigma \lambda \end{array} \right)
\left(\begin{array}{cc} 0 &0 \\ -\sigma^2 \gamma & -\sigma- \sigma\eta'\sigma -\sigma^2 \eta  \end{array}\right)
\left( \begin{array}{cc} 1 \\ -\sigma \lambda \gamma \end{array} \right)\\
&=& \left( \begin{array}{cc} \beta\sigma\lambda\sigma^2 \gamma  &\beta \sigma\lambda(\sigma+ \sigma\eta'\sigma +\sigma^2 \eta   ) \end{array} \right)
\left( \begin{array}{cc} 1 \\ -\sigma \lambda \gamma \end{array} \right)\\
&=& \beta\sigma\lambda\underline{\sigma^2} \gamma - \beta \sigma\lambda\underline{(\sigma+ \sigma\eta'\sigma +\sigma^2 \eta   )\sigma \lambda} \gamma\\
&=& \beta \sigma\lambda (\sigma^2 - \sigma^2\lambda - \sigma \eta' \sigma^2\lambda - \sigma^2 \underline{\eta} \sigma \lambda) \gamma\\
&=&  \beta \sigma\lambda (\underline{\sigma^2 - \sigma^2\lambda} - \sigma \eta' \sigma^2\lambda - \sigma^2\eta' \sigma \lambda \underline{- \sigma^2\eta'' \sigma \lambda}) \gamma\\
&=&  \beta \sigma\lambda ( \sigma^2\underline{(1-\lambda - \eta''\sigma \lambda)}- \sigma \eta' \sigma^2\lambda - \sigma^2\eta' \sigma \lambda ) \gamma\\
&\overset{(*)}{=}&  \beta \sigma\lambda ( - \sigma \underline{\eta' \sigma}\sigma \lambda - \sigma\underline{\sigma\eta'} \sigma \lambda ) \gamma\\
&=& -\beta \sigma\lambda \sigma^2\lambda \gamma.
  \end{array}$$
  established a homotopy between $gf$ and $1$, where equation $(*)$ use the fact that $$1-\lambda - \eta''\sigma\lambda=1-\lambda (1+\eta''\sigma)=0.$$

  This shows (a), (c) and part of (b).

\medskip

  Let us prove the remaining part of (b).  Obviously $f$ is injective and $g$ is surjective.

Now we show that the maps
$$  \left(\begin{array}{c} \beta \sigma \\ 1+\eta''\sigma \end{array}\right): (C_*',\eta')\rightarrow \mathrm{Ker}(g)  \ \ \mathrm{and}\ \
(\sigma\gamma  \  \  1+\sigma\eta''): \mathrm{Cok}(f)\rightarrow(C_*', \eta')
$$
are isomorphisms of complexes with inverse maps
$$\begin{array}{cclcl} \left(\begin{array}{c} \beta \sigma \\ 1+\eta''\sigma \end{array}\right)^{-1}&=& (0 \ \ \lambda) &:&\mathrm{Ker}(g)\rightarrow (C_*',\eta'),\\
(\sigma \gamma  \  \  1+\sigma\eta'')^{-1}&=&
\left(\begin{array}{c} 0 \\ \overline\lambda \end{array}\right)&:&(C_*', \eta')\rightarrow \mathrm{Cok}(f)
\end{array}
$$
respectively,
where $\overline\lambda=1-\sigma\lambda\eta''$ is the inverse of $1+\sigma\eta''$

 In fact,
$$
g\circ \left(\begin{array}{c} \beta \sigma \\ 1+\eta''\sigma \end{array}\right) =\beta \sigma-\beta \sigma\lambda (1+  \eta'' \sigma)=0,$$
 so  $$\left(\begin{array}{c} \beta \sigma \\ 1+\eta''\sigma \end{array}\right):  (C_*',\eta')\rightarrow Ker(g)$$ is well-defined.
We have
$$\begin{array}{rcl}
  d\circ \left(\begin{array}{c} \beta \sigma \\ 1+\eta''\sigma \end{array}\right)&=&\left(\begin{array}{cc} \alpha&\beta \\ \gamma&\eta \end{array}\right)
\left(\begin{array}{c} \beta \sigma \\ 1+\eta''\sigma \end{array}\right)\\
&=&\left(\begin{array}{c}\alpha\beta \sigma+\beta(1+\eta''\sigma) \\ \gamma \beta \sigma+\eta(1+\eta''\sigma) \end{array}\right)
=\left(\begin{array}{c} -\beta \eta \sigma+\beta(1+\eta''\sigma) \\ -\eta^2 \sigma+\eta(1+\eta''\sigma) \end{array}\right)\\
&=&\left(\begin{array}{c} \beta(1+\eta''\sigma-\eta \sigma) \\ \eta(1+\eta''\sigma-\eta \sigma) \end{array}\right)
=\left(\begin{array}{c} \beta \sigma\eta' \\ \eta \sigma\eta' \end{array}\right)
 = \left(\begin{array}{c} \beta \sigma\eta' \\ \eta'+\eta''\sigma\eta' \end{array}\right)\\
 &=&\left( \begin{array}{c} \beta \sigma\\ 1+\eta'' \sigma \end{array}\right)\eta'.
\end{array}
$$
Hence it is a chain map.
Now  $$(0 \ \ \lambda)
 \left( \begin{array}{c} \beta \sigma \\ 1+\eta'' \sigma \end{array}\right)=\lambda(1+\eta''\sigma)=1;$$
 for the other equality, observe that an  arbitrary  element of $\mathrm{Ker}(g)$ has the form
 $\left( \begin{array}{c} \beta \sigma \lambda(y) \\ y\end{array}\right)$ with $y\in C'_*$ and so
 $$\left( \begin{array}{c} \beta \sigma \\ 1+\eta'' \sigma \end{array}\right)(0 \ \ \lambda)\left( \begin{array}{c} \beta \sigma \lambda(y) \\ y\end{array}\right)
=\left( \begin{array}{c} \beta \sigma \lambda(y)\\ (1+\eta'' \sigma )\lambda(y) \end{array}\right)
=\left( \begin{array}{c} \beta \sigma\lambda(y) \\ y \end{array}\right).
$$
We have shown that $ \left( \begin{array}{c} \beta \sigma \\ 1+\eta'' \sigma \end{array}\right)$ is an isomorphism of chain complexes.

\medskip

For the   statements about $(\sigma\gamma  \  \  1+\sigma\eta'') :  \mathrm{Cok}(f)\rightarrow(C_*', \eta')$,
we have $$  (\sigma \gamma \ \ 1+\sigma\eta'')
\left( \begin{array}{c} 1 \\ -\sigma\lambda \gamma  \end{array}\right)
 = \sigma \gamma -(1+\sigma\eta'') \sigma\lambda \gamma
 = \sigma \gamma-\sigma\lambda\gamma-\sigma(1-\lambda)\gamma \\
 = 0.
 $$
Hence $(\sigma\gamma  \  \  1+\sigma\eta'') :  \mathrm{Cok}(f)\rightarrow(C_*', \eta')$ is well-defined.

As $$(\sigma \gamma \ \ 1+\sigma\eta'') \left( \begin{array}{cc} \alpha&\beta\\ \gamma &\eta \end{array} \right)
=(\sigma \gamma \alpha+(1+\sigma\eta'') \gamma \ \ \ \sigma \gamma \beta+(1+\sigma\eta'') \eta )=( \eta'\sigma\gamma \ \ \eta' \sigma\eta)$$
and
 $\eta' (\sigma\gamma \ \ 1+\sigma\eta'')=(\eta'\sigma\gamma \ \ \eta'+\eta'\sigma\eta''),$
together
with the fact that $$\eta'\sigma\eta=\eta-\sigma\eta'\eta=\eta-\sigma\eta'\eta''=\eta'+(1-\sigma\eta')\eta''
=\eta'+\eta'\sigma\eta'',$$
we see that  $(\sigma\gamma \ \ 1+\sigma\eta'')$ is a chain map.

We have  $$(\sigma\gamma  \ \ 1+\sigma\eta'')\left(  \begin{array}{c} 0\\ \overline\lambda \end{array} \right)=(1+\sigma\eta'')\overline\lambda=1$$
and for $x\in C_*, y\in C_*'$,
$$\left(  \begin{array}{c} 0\\ \overline\lambda \end{array} \right)(\sigma\gamma \ \ 1+\sigma\eta'') \overline{\left(  \begin{array}{c} x \\ y \end{array}  \right)}=\overline{\left(  \begin{array}{c} 0\\ \overline\lambda \sigma\gamma(x)+y \end{array} \right)}=\overline{\left(  \begin{array}{c} x\\ y\end{array} \right)}$$
where $\overline{ \left(  \begin{array}{c} x \\ y \end{array}  \right) } \in Cok(f)$.
Then  $(\sigma\gamma  \  \  1+\sigma\eta'') :  \mathrm{Cok}(f)\rightarrow(C_*', \eta')$  is an isomorphism of complexes.
\end{Proof}

\begin{Rem}\label{Rem: h2=0}

  In   Theorem~\ref{Thm: reformulation-algebraic-morse}, if $\sigma^2=0$, then $h^2=0$ and the original datum is a SDR datum, so the perturbed datum is also an SDR datum by Theorem~\ref{Thm: perturbation lemma}. In this case, we have  $gf=1$ and the two short exact sequences of Theorem~\ref{Thm: reformulation-algebraic-morse} (b) split each other.

 \end{Rem}

The following corollary (the so-called \textit{killing-contractible-complexes lemma}) is well-known in homological algebra and we cite it from the textbook \cite[Lemma 2.1.6]{Loday}. Note that this result is stated over a commutative ring $k$ in \cite{Loday} but it is obviously true over any associative ring. The form we present here is more precise than \cite[Lemma 2.1.6]{Loday}, where it is only stated as a quasi-isomorphism.

\begin{Cor}  \label{killing-contractible-comlexes} Let $R$ be an associative ring, and let $$ \cdots\longrightarrow
C_{n}\oplus C_{n}'\stackrel{d=\left(
\begin{array}{cc} \alpha&\beta \\ \gamma&\eta \end{array}\right)}{\longrightarrow}
C_{n-1}\oplus C_{n-1}'\longrightarrow\cdots$$ be a chain complex of $R$-modules such that $(C_*',\eta)$ is a complex which is contractible with contracting homotopy $\sigma: C_n'\longrightarrow C_{n+1}'$. Then the following inclusion of complexes is a homotopy equivalence: $$\left(
\begin{array}{c} 1 \\ -\sigma\gamma \end{array}\right): (C_*,\alpha-\beta \sigma\gamma)\hookrightarrow (C_*\oplus C_*',d).$$
\end{Cor}

\bigskip

\section{Algebraic Morse theory}\label{Sect: morse theory}

In this section, we present a version of
  algebraic Morse theory which is more  general than existing ones in the literature.

Let $R$ be a ring, all modules will be (left) $R$-modules.

 Let $X_*$ be the following complex of $R$-modules:
$$\cdots\to X_{n+1}\stackrel{d_{n+1}}{\to} X_n \stackrel{d_{n}}{\to} X_{n-1}\to\cdots$$
Suppose that for each $n\in \Z$,  there exists a decomposition into direct sums of submodules
$$X_n=\oplus_{i\in I_n}X_{n, i}.$$
So $d_n:X_n\to X_{n-1}$ has a matrix presentation  $d_n=(d_{n, ji})$ with $i\in I_n, j\in I_{n-1}$ and
where $d_{n, ji}: X_{n, i}\to X_{n-1, j}$ is a homomorphism of modules.

We shall construct a  weighted  quiver  $Q=Q_{X_*}$  as follows:
\begin{itemize}

\item[(Q1)]
The  vertices are the pairs   $(n, i)$ with $ n\in \Z, i\in I_n$;

\item[(Q2)] if a map $d_{n, ji}$ with $i\in I_n, j\in I_{n-1}$ does not vanish, then draw an arrow from $(n, i)$ to $(n-1, j)$;

\item[(Q3)] for an arrow in (Q2), its  weight  is just the map  $d_{n, ji}$.
 \end{itemize}

 A \textit{partial  matching} is a full subquiver $\M$ of $Q$ such that
 \begin{itemize}
 \item[(M1)]each vertex in $Q$  belongs to at most  one arrow of $\M$;

 \item[(M2)] each arrow in $\M$  has its weight  invertible as a homomorphism of modules.
 \end{itemize}

 Given a  partial matching  $\M$, we can construct a new weighted quiver $Q^\M$ with additional dotted arrows as follows:
 \begin{itemize}
 \item[(QM1)]Keep everything for all  arrows  which are not in $\M$ (they will be called \textit{thick arrows});

 \item[(QM2)] For an  arrow in $\M$, replace it by  a new \textit{dotted arrow} in the reverse direction and the weight of this new arrow is the negative  inverse of the weight of the original arrow.
 \end{itemize}

Note that both $Q$ and $Q^\M$   have no multiple arrows.

We need to fix some notations:
For $n\in \mathbb{Z}$,  denote   $\V_n=\{(n, i): i\in I_n\}$.  Call a vertex of $Q$   a \textit{critical vertex} (with respect to $\M$), if it is not incident to  any arrow in $\M$. Denote by $\V_n^\M$ the set of critical vertices in $\V_n$. Denote by $\U_n\subseteq \V_n$ the set of vertices which appear as the starting vertex of an arrow in $\M$,  by $\D_n\subseteq \V_n$ the set of vertices which appear as the ending  vertex of an arrow in $\M$.  Clearly, we have $\V_n=\V_n^\M\bigcup \U_n\bigcup \D_n$ (disjoint union).
 For a path $p$ in $Q^\M$, we can define a composition which is the composition of all maps appearing as the weights of all arrows in $p$, denote this composition as $\varphi_p^\M$.

 A path in $Q^\M$ is called zigzag if dotted arrows and  thick arrows  appear  alternately.  For any two vertices $(n, i)\in \V_n$ and $(m, j)\in \V_m$, denote by    $\Path^\M((n, i), (m, j))$ the set of all zigzag paths from $(n, i)$ to $(m, j)$  in $Q^\M$; notice that in this case, we have necessarily $m=n, n-1$ or $n+1$.
For  any two vertices $(n, i)\in \D_n$ and $(n, j)\in \V_n$, denote by    $\Path^\M_1((n, i), (n, j))$ the set of all zigzag paths from $(n, i)$ to $(n, j)$  in $Q^\M$   which  begin  with  a dotted arrow  (and which necessarily  end   with  a thick arrow).  Similarly, we denote by $\Path^\M_2((n, i), (n, j))$ the set of all zigzag paths from $(n, i)$ to $(n, j)$  in $Q^\M$   which  begin  with  a thick arrow  (and which necessarily  end   with  a dotted arrow).

 We impose a \textit{local  finiteness hypothesis} (LFH)  which is more general than conditions imposed previously; for the comparison between this condition with those of \cite{Kozlov,S06,JW}, see Remarks~\ref{Remark: comparison Kozlov}, \ref{Remark: comparison S} and \ref{Remark: comparison JW}.  A \textit{Morse matching} is a partial  matching which satisfies the following local finiteness hypothesis:
 \begin{itemize}

 \item[(LFH)] Given  an arbitrary vertex  $(n, i)\in \D_n$, for each vertex  $(n, j)\in  \V_{n}$ and
     for each element $x\in X_{n, i}$, the sum
     $$\sum_{ p\in  \Path^\M_1 ((n, i),  (n, j) ) } \varphi_p^\M(x)$$
     exists (for instance, it may be a finite sum or it is convergent in a certain norm); moreover,
     the number of vertices $(n, j)\in \V_{n}$ such that  $$ \sum_{p\in  \Path^\M_1 ((n, i),  (n, j) )} \varphi^\M_p(x) \neq 0$$  is finite.

 \end{itemize}

 \begin{Rem}  Since for each element $x\in X_{n, i}$ the value $d_{n}(x)$ is a finite sum in $X_{n-1}=\oplus_{j\in I_{n-1}}X_{n-1, j}$, it  is easy to see that the condition (LFH) will imply the following three conditions:

 \begin{itemize}
 \item[(LFH1)]
Given an arbitrary vertex  $(n, i)\in \V_n^\M$, for each vertex  $(n-1, j)\in \V_{n-1}^\M$  and
     for each element $x\in X_{n, i}$, the sum
     $$\sum_{ p\in \Path^\M((n, i), (n-1, j)) } \varphi_p^\M(x)$$
     exists; moreover,
     the number of vertices $(n-1, j)\in \V_{n-1}^\M$ such that  $$ \sum_{p\in \Path^\M((n, i), (n-1, j))} \varphi^\M_p(x)  \neq 0$$  is finite.

 \item[(LFH2)]
 Given an arbitrary vertex  $(n, i)\in \D_n$, for each vertex  $(n+1, j)\in \U_{n+1}$  and
     for each element $x\in X_{n, i}$, the sum
     $$\sum_{ p\in \Path^\M((n, i), (n+1, j)) } \varphi_p^\M(x)$$
     exists; moreover,
     the number of vertices $(n+1, j)\in  \U_{n+1} $ such that  $$ \sum_{p\in \Path^\M((n, i), (n+1, j))} \varphi^\M_p(x) \neq 0$$  is finite.

\item[(LFH3)]  Given an arbitrary vertex  $(n, i)\in \V_n^\M$, for each vertex  $(n, j)\in \V_{n}$, and
     for each element $x\in X_{n, i}$, the sum
     $$\sum_{ p\in    \Path^\M_2((n, i), (n, j)) }   \varphi_p^\M(x)$$
     exists; moreover,
     the number of vertices $(n, j)\in \V_{n}$ such that  $$ \sum_{p\in    \Path^\M_2 ((n, i), (n, j))}   \varphi^\M_p(x) \neq 0$$  is finite.

 \end{itemize}

 \end{Rem}

The following proposition is more readily applicable to check Morse matching, which is a reformulation of Sk\"{o}ldberg's sufficient condition in Remark~\ref{Remark: comparison S}. 

\begin{Prop}\label{Prop: sufficient condition for Morse matching}
 Let $\M$ be a partial matching of $Q$. If any zigzag path from $(n,i)$ is of finite length for each vertex $(n,i)$ in $ Q^{\M},$ then $\M$ is a Morse matching.
\end{Prop}

\begin{Proof}
 Let $(n,i) \in \D_n$. We first prove that the set of all zigzag paths from $(n,i)$ which begin with a dotted arrow is a finite set. Assume that there are infinite such zigzag paths in $Q^{\M}$. The condition (M1) of the partial matching $\M$ guarantees that the first dotted arrows of these zigzag paths all coincide with an arrow of the form $(n,i) \dashrightarrow (n+1,j)$.
Since each term $X_n$ in the complex $X_*$ is a direct sum, there are only finitely many thick arrows leaving from $(n+1,j)$. So there exists at least one thick arrow $(n+1,j) \to (n,k)$ such that there are infinite zigzag paths which begin with
$$(n,i) \dashrightarrow (n+1,j)\to (n,k).$$
Repeating the above process for the zigzag paths from $(n,k)$ and by induction, we obtain a zigzag path from $(n,i)$ of infinite length, a contradiction. Hence the set
of all zigzag paths  from $(n,i)$ which begin with a dotted arrow is a finite set. In particular, for each element $x\in X_{n, i}$, the set
$$\{p\mid p\text{ is a zigzag path which begins with a dotted arrow from (n,i)}, \varphi^\M_p(x)\neq 0\}$$
 is a finite set. Therefore the condition (LFH) follows.

\end{Proof}

 Given a Morse matching $\M$, we can construct a new complex $( {X}_*^\M, d^\M)$ as follows:

The complex ${X}_*^\M$ has its $n$-th component $X_n^\M=\oplus_{(n, i)\in \V_n^\M} X_{n, i}$ and the differential
$d_n^\M: X_n^\M\to X_{n-1}^\M$ has the matrix presentation
$d_n^\M=(d_{n, ji}^\M)$ with $(n, i)\in \V_n^\M, (n-1, j)\in \V_{n-1}^\M$ and where
$d_{n, ji}^\M: X_{n, i}\to X_{n-1, j} $  is defined to be
$$d_{n, ji}^\M=\sum_{p\in  \Path^\M((n, i), (n-1, j))} \varphi^\M_p.$$
Note that $d^\M$ exists by (LFH1).

Now we can state the main result of algebraic Morse theory, which  contains as special cases and refines  all versions appearing in the literature \cite{Kozlov,S05,JW}.

\begin{Thm}\label{Thm: main result of algebraic Morse theory}
\begin{itemize}

\item[(a)] Within the above setup,
  $({X}_*^\M, d^\M)$ is a complex.

\item[(b)]   Define maps \begin{align*}
f_n: X_n^\M & \rightarrow X_n \\
 x \in X_{n,i} & \mapsto f_n(x):=x+ \sum_{(n, j)\in \U_n}\sum_{p\in  \Path^\M_2   ((n,i), (n, j))} \varphi^\M_p(x),
\end{align*}and

\begin{align*}
g_n: X_n & \rightarrow X^\M_n \\
 x \in X_{n,i} & \mapsto g_n(x):=\left\{\begin{array}{ll} \sum_{(n, j)\in \V^\M_n}\sum_{p\in   \Path^\M_1   ((n,i), (n, j)) } \varphi^\M_p(x),& (n, i)\in \D_n\\
  x,  &  (n,i) \in \V_n^\M \\
 0 & (n, i)\in \U_n\end{array}\right.
\end{align*}
 Then $f_*: {X}_*^\M   \rightarrow X_* $ and $ g_*: X_*   \rightarrow {X}_*^\M $ are chain maps which are homotopy equivalent:
$gf=\Id_{{X}_*^\M}$ and $fg\sim \Id_{X_*}$ via the homotopy
\begin{align*}
{\theta_n}: X_n & \rightarrow X_{n+1} \\
  x \in X_{n,i} & \mapsto {\theta_n}(x):= \left\{\begin{array}{ll} \sum_{(n+1, j)\in \U_{n+1}}\sum_{p\in \Path^\M((n,i), (n+1, j)) } \varphi^\M_p(x),& (n, i)\in \D_n\\
   0 & otherwise\end{array}\right.
\end{align*}

    \item[(c)] We have a decomposition $$(X_*,d) \cong  (X_* ^\M, d^\M)  \oplus (Y_*,d^Y)$$ where $(Y_*,d^Y)$ is a null homotopic  complex.
\end{itemize}
\end{Thm}

\begin{Proof}  We will use Theorem~\ref{Thm: reformulation-algebraic-morse} to prove the result.

Construct a chain complex of $R$-modules $({C_{n}\oplus C_{n}'},d)$, $d=\left(  \begin{array}{cc}\alpha & \beta \\ \gamma &\eta
\end{array}\right)$,
as follows:

 let
$${C_n}=\oplus_{(n, i)\in \V_n^\M} X_{n, i}, {C_n'}=\oplus_{(n, i)\in \V_n\backslash \V_n^\M} X_{n, i},$$ $\alpha,\beta,\gamma,\eta$ be the maps defined by corresponding arrows in $Q$ respectively, $\eta'$ be the map defined by arrows in $\M$, $\eta''=\eta-\eta'$ be the weight of the arrows  which are not in $\M$;
for an arrow $d_{n, ji}: X_{n, i}\to X_{n-1, j}$ lying in $\M$, then  ${\sigma}=d_{n, ji}^{-1}:X_{n-1, j}\rightarrow X_{n, i}$ is the homotopy for the ``piece" $d_{n, ji}: X_{n, i}\to X_{n-1, j}$, so the weight of the dotted arrow is $-{\sigma}$.
Note that we have always ${\sigma}^2=0$.

It's obvious that $({C_{n}\oplus C_{n}'},d)$ is just the original  chain complex $X_*$; the fact that $\eta'^2 =0$ and that ${\sigma}$ is contracting homotopy of $({C_n'}, \eta')$ follows from the definition of Morse matching; it is interesting to see that powers of $\eta''{\sigma}$ are zigzag paths in $Q^\M$, so for each $x\in X_{n, i}$ with $(n, i)\in \D_n$,  $$\lambda(x):=(1+\eta''{\sigma})^{-1}(x)=\sum\limits_{i=0}^{\infty} (\eta''(-{\sigma}))^i(x)$$ exists  by  the axiom (LFH).

 Note that the existence of  $f$, $g$ and ${\theta}$ follows from (LFH3), (LFH) and (LFH2) respectively.  It is  not difficult to check that $f$ and $g$ are exactly the chain maps $f$ and $g$ in Theorem~\ref{Thm: reformulation-algebraic-morse}, and that  the homotopy ${\theta}$ is just $${h_{\infty}=\left( \begin{array}{cc} 0 &0 \\ 0 & -\sigma \lambda \end{array} \right)}$$ in Theorem~\ref{Thm: reformulation-algebraic-morse}. Now by applying Theorem~\ref{Thm: reformulation-algebraic-morse} and Remark~\ref{Rem: h2=0}, and by noting that the differential in the new complex ${X}_*^\M$ is just $\alpha-\beta{\sigma}\lambda\gamma$, we get the desired  result.

\end{Proof}

\begin{Rem} By Remark~\ref{Rem: h2=0}, the above theorem in fact gives an {SDR datum}
$$ \xymatrix{ {X}_*^\M\ar@<1ex>[r]^{f} & X_* \ar@<1ex>[l]^{g} \ar@(ur,dr)^{{\theta}}}. $$

\end{Rem}

\begin{Rem}\label{Remark: comparison Kozlov} In the  paper  \cite{Kozlov}, D.~N.~Kozlov defined a Morse matching to be a partial matching satisfying the following conditions:
\begin{itemize}
\item[(K1)]  Each $I_n$ is finite;

\item[(K2)]  for each vertex $(n, i)\in \D_n$,  $\Path^\M_1 ((n, i), (n, i))=\emptyset$, i.e. there is no zigzag path, which begins with a dotted arrow and ends with a thick arrow,  from  $(n, i)$ to itself.
\end{itemize}

In fact, D.~N.~Kozlov also asked that $R$ is commutative ring and that $X_n=0$ for $n\ll 0$. 
It is obvious that the first condition is unnecessary and the latter one is superfluous by \cite[Remark 2]{Kozlov}.

It is easy to see that whenever the conditions (K1)(K2) hold, the condition in Proposition~\ref{Prop: sufficient condition for Morse matching} is satisfied, and hence the    axiom  (LFH) holds.  
On the other hand, the condition of Proposition~\ref{Prop: sufficient condition for Morse matching}, together with the finiteness of each $I_n$, implies conditions (K1)(K2).
\end{Rem}

\begin{Rem}\label{Remark: comparison S}  In the article \cite{S06},   E.~Sk\"{o}ldberg defined a Morse matching to be a partial matching satisfying the following condition:
\begin{itemize}
	\item[(S)] There is a well-founded partial order  $\prec$ on  each $I_n$  such that  $c\prec a $  whenever there is a zigzag path $a\to b \to c$ in $Q^{\M}$ (note that one of the arrows should be dotted). 

\end{itemize}
 Here we note that the condition (S) is equivalnet to the  condition in Proposition~\ref{Prop: sufficient condition for Morse matching}, and hence implies the  axiom  (LFH).   Indeed,   each zigzag path  
  	$$  a_1 \to b_1\to a_2 \to b_2 \to \cdots \to a_n \to b_n \to \cdots     $$
  	gives rise to a descending sequence $a_1 \succ a_2 \succ   \cdots$ with respect to $\prec$. Therefore, the condition (S) implies the condition in Proposition~\ref{Prop: sufficient condition for Morse matching}.  
  	On the other hand, given a descending sequence  $a_1 \succ a_2 \succ   \cdots$ where all $a_i$ belong to $I_n$, we may assume without loss of generality that it corresponds to  a  path
  	$$   \xymatrix@R=1ex{ & a_1'=a_1 \ar@{.>}[ld]& \\
  	c_1 \ar[rd]& &  \\
  	& a_2' & \\
  	& \vdots  & \\
  	& a_{m-1}'\ar@{.>}[ld] &\\
  c_{m-1}	\ar[rd]& & \\
  	& a_m'\ar[rd]& \\
  	& & d_m \ar@{.>}[ld]\\
  	&a_{m+1}'& \\
  	& \vdots &  }  $$
   in $Q^{\M},$ which contains a zigzag subpath starting at $a_m'.$ Hence, the condition in  Proposition~\ref{Prop: sufficient condition for Morse matching} implies the condition (S).

\end{Rem}

\begin{Rem}\label{Remark: comparison JW}  In the article \cite{JW}, M.~J\"{o}llenbeck and  V.~Welker  defined a Morse matching to be a partial matching satisfying the following conditions:
\begin{itemize}
\item[(JW1)] $Q^\M$ contains no directed cycles;

\item[(JW2)]  the finiteness condition in \cite[Definition 2.3]{JW} holds.

\end{itemize}

It is not hard to see that whenever the conditions  (JW1)(JW2) hold, the condition in Proposition~\ref{Prop: sufficient condition for Morse matching} is satisfied, and hence the    axiom  (LFH) holds.
\end{Rem}

\bigskip

\section{Two-sided Anick resolutions}\label{Sect: Two-sided Anick resolutions}

In this section, we generalise the construction of two-sided Anick resolutions of E.~Sk\"{o}ldberg in \cite{S06} from one vertex algebras to algebras given by quotients of path algebras of quivers.

\medskip



Let $k$ be a fixed field.  Let $Q=(Q_0,Q_1,s,t)$ be a finite quiver with vertex set $Q_0$ and  arrow set $Q_1$, where $s: Q_1\to Q_0$ (resp. $t: Q_1\to Q_0$) gives the starting vertex (resp. the target vertex) of an arrow.
    We will write paths from left to right, that is, the notation $p=\alpha_1\alpha_2\cdots\alpha_r$ means that for $1\leq i\leq r-1$, $s(\alpha_{i+1})=t(\alpha_i)$. The length $l(p)$ of the above path $p$ is defined to be $r$ and the vertices are viewed as paths of length $0$.   For $n\geq 0$, $Q_n$ denotes the set of all paths of length $n$ and $Q_{\geq n}$  is  the set of all paths with   length at least $n$. Denote by $kQ$ its path algebra, that is the space generated by all paths of finite length and endowed with
the multiplication given by concatenation of paths.

 Let us first briefly recall the Gr\"{o}bner-Shirshov basis theory for a path algebra $kQ$  following E.~Green's paper \cite{Green}. Let $\mathcal{B}:=Q_{\geq 0}$ denote the set of all  finite (directed) paths in $Q$. Then $\mathcal{B}$ is a multiplicative $k$-basis of $kQ$.  Write $\mathcal{B}_+=\mathcal{B}\setminus Q_0$.  Fix an admissible well-order  $\prec$  on $\mathcal{B}$, that is, a well-order on $\mathcal{B}$ which is compatible with multiplication. For instance, we can take a left length-lexicographic order extending some total order on the arrows. For a linear  combination $r$  of paths, its tip $\mathrm{Tip}(r)$ is by definition the maximal monomial appearing with nonzero coefficients in $r$. For a nonempty  subset $X$ of $kQ$, put
 $\mathrm{Tip}(X)=\{\mathrm{Tip}(r)\ |\ r\in X, r\neq 0\}$.

 Let $I$ be a two-sided ideal in $kQ$ contained in $kQ_{\geq 2}$.
    Write  $\NonTip (I)$ for the complement set of ${\Tip} (I)$ in $\mathcal{B}$. Then there exists a decomposition of vector spaces
 $$kQ=I\oplus \text{Span}_{k}(\NonTip(I)).$$
 So $\NonTip(I)$ is a basis of the quotient algebra $A=kQ/I$. Recall that a Gr\"{o}bner-Shirshov basis of $I$ with respect to  the admissible order  $\prec$ is a   subset $\mathcal{G}\subseteq I$ such that $W:=\Tip(\mathcal{G})$ generates the initial ideal $\langle {\Tip} (I)\rangle$. Note that  in this case  $I=\langle \mathcal{G}\rangle$.  A Gr\"{o}bner-Shirshov basis $\mathcal{G}$ for the ideal $I$ is reduced if the following three conditions are satisfied: \begin{itemize}
\item[(R1)] For any $g\in \mathcal{G}$, the coefficient of $\Tip(g)$ is $1$;
\item[(R2)] For any $g\in \mathcal{G}$, $g-\Tip(g)\in \text{Span}_{k}(\NonTip(I))$;
\item[(R3)] No element in $W=\Tip(\mathcal{G})$ is a factor of another element in $W$.
\end{itemize}
\noindent It is easy to see that under the given admissible order, $I$ has a unique reduced Gr\"{o}bner-Shirshov basis, and in this case $W$ is a minimal generator set of $\langle {\Tip} (I)\rangle$; moreover, $b\in \mathcal{B}$ lies in $\NonTip(I)$ if and only if $b$ is not divided by any element of $W$. In the following, we always assume that $W=\Tip(\mathcal{G})$ for a reduced Gr\"{o}bner-Shirshov basis $\mathcal{G}$ of $I$.

Similar as in \cite{S06}, we define a new quiver $Q_W=(V,E)$ (with respect to $W$), called  the Ufnarovski\u\i~graph \cite{U},  with vertex set $V$ and arrow set $E$ as follows:
 $$V=Q_0\cup Q_1 \cup \{u\in \mathcal{B}~|~u\text{ is a proper right factor of some }v\in W \},$$
 and $E$ is the union of
$ \{ e\to x\ |\ e\in Q_0, x=ex\in Q_1 \}$
 with $$\{u\to v\  |\  uv\in \mathcal{B}, uv\in \langle \Tip(I)\rangle, w\notin \langle \Tip(I)\rangle\text{ for all proper left factors }w\text{ of }uv\}.$$
Note that the above condition for the arrow $u\to v$ is equivalent to the following: $uv\in \mathcal{B}$ and $uv$ has a unique factor $w\in W$ which is a right factor of $uv$. Clearly the above condition implies the following: if $u\to v$ is an arrow in $Q_W$, then $u\to v_1$ can not be an arrow in $Q_W$ for any proper left factor $v_1$ of $v$.  The set of $i$-chains  $W^{(i)}, i\ge 0$ (also called Anick chains)  consists of all sequences $(v_1, \cdots ,v_i,v_{i+1})$ in $\mathcal{B}_+^{i+1}$ such that
 $$e\rightarrow v_1\rightarrow \cdots \rightarrow v_i\rightarrow v_{i+1}$$
 is a path in $Q_W$, where $e\in Q_0$ and $\mathcal{B}_+=\mathcal{B}\setminus Q_0$; by convention $W^{(-1)}=Q_0$.  For $w=(w_1,\cdots, w_m)\in W^{(m-1)}$,   the length of the path $w_1\cdots w_m$ is  called  the degree of   $w$ and $m$ is called the weight of $w$.  In order to have an intuition, let us consider a concrete example of the Ufnarovski\u\i~graph.

\begin{exam} \label{exam:Uf-graph} Let $A=kQ/I$ be a finite dimensional algebra defined by the following quiver \unitlength=1.00mm
\special{em:linewidth 0.4pt} \linethickness{0.4pt}
\begin{center}
\begin{picture}(12.00,14.00)(11,0)
\put(5,5){$1$} \put(20,5){$2$} \put(35,5){$3$}
\put(-5,5){$Q:$} \put(9,7){\vector(1,0){8.0}}
\put(17,5){\vector(-1,0){8.0}} \put(24,7){\vector(1,0){8.0}}
\put(32,5){\vector(-1,0){8.0}}\put(12,8){$a$}
\put(12,1){$a'$} \put(27,8){$b$}\put(27,1){$b'$}
\end{picture}
\end{center}with relations
$$ab=b'a'=a'a-bb'=0.$$
Consider the left length-lexicographic order with $a>b>b'>a'$. Then $\mathcal{G}=\{ab, b'a', bb'-a'a, aa'a, a'aa'\}$ is a reduced Gr\"{o}bner-Shirshov basis of $I$ with respect to this order. Therefore we have the following: $$W=\Tip(\mathcal{G})=\{ab, b'a', bb', aa'a, a'aa'\}, Q_0=\{1,2,3\}:=\{e_1,e_2,e_3\}, Q_1=\{a, b, a', b'\},$$
$$\{u \text{ is a proper right factor for some }v\in W\}\setminus Q_1=\{a'a,aa'\}, \NonTip(I)\setminus Q_0=\{a,b,b',a',a'a,aa',b'b\}.$$
The associated Ufnarovski\u\i~graph $Q_W$ is
$$
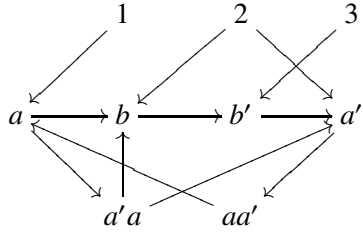
\xymatrix{& 1\ar[ld] & 2\ar[ld]\ar[rd] & 3\ar[ld] \\
 a\ar[r]\ar[rd] &b\ar[r] & b'\ar[r] & a'\ar[ld]\\
 & a'a\ar[u]\ar[rru] & aa'\ar[llu] & }$$
\end{exam}

Now let $A=kQ/I$ be  an arbitrary algebra as before.  Then $E=\oplus_{e\in Q_0}ke$ is a semisimple subalgebra of $A$ such that $A=E\oplus A_+$ as spaces, where $A_+=\text{Span}_{k}\{\NonTip(I)\setminus Q_0\}$. Recall the reduced two-sided bar resolution $B(A,A)$ of the algebra $A$ in the sense of C.~Cibils \cite{Cibils}:
 $$B(A,A)_0=A\otimes_EA, \text{ and for }n\geq 1, B(A,A)_n=A\otimes_E {(A_+)}^{\otimes_E n}\otimes_E A \cong A^e \otimes_{E^e} {(A_+)}^{\otimes_E n},$$
 where $A^e=A\otimes_kA^{op}$ and  $E^e=E\otimes_kE^{op}$, and the differential is defined by (for $n\geq 1$)
 \begin{equation}\label{Eq: differential in reduced bar resolution}d([a_1 | \cdots |a_n])=a_1 [a_2 |\cdots| a_n]+\sum\limits_{i=1}^{n-1} (-1)^i [a_1| \cdots|a_i a_{i+1}|\cdots |a_n]+(-1)^n [a_1| \cdots | a_{n-1}]a_n,\end{equation}
 where $[a_1|\cdots |a_n]=1\otimes a_1\otimes \cdots \otimes a_n\otimes 1$.
  We decompose the above resolution as follows:
 $$B(A,A)_n=\bigoplus A^e \otimes_{E^e} k[w_1 | \cdots | w_n]=\bigoplus A^e \cdot [w_1 | \cdots | w_n],$$
 where the sum is taken over all the sequences $(w_1,\cdots ,w_n)$ such that all $w_i\in \NonTip(I)\setminus Q_0$ and $w_1\cdots w_n$ is a path in $Q$.  Note that if $w_1\cdots w_n=e_iw_1\cdots w_ne_j$, then $A^e \cdot [w_1 | \cdots | w_n]$ is isomorphic to the indecomposable projective $A^e$-module $Ae_i\otimes_ke_jA$.  We write $(w_1,\cdots ,w_n)$ instead of $[w_1 | \cdots | w_n]$ for the vertices in the decorated quiver $Q_B=Q_{B(A,A)}$ (cf. Section 3).

  In general, an arrow in $Q_B$ may contain the information of several terms in the expression (\ref{Eq: differential in reduced bar resolution}) of the differential (for example $d( [a])= a [ ] 1 -1 [ ] a $ and its corresponding arrow is $(a)\overset{a\otimes 1 -1 \otimes a}{\to} e_1 \in Q_{0} $). For the sake of clarity, it is necessary to view each term of the expression of the differential $d$ as an arrow. That is to say we will use a new weighted quiver $\overline{Q_{B}}$ to construct the two-sided Anick resolution of $A$ via Theorem~\ref{Thm: main result of algebraic Morse theory} instead of $Q_B.$ The weighted quiver $\overline{Q_B}$ has the same vertex set as $Q_B$. And the arrows  of $\overline{Q_{B}}$ are listed as follows. For a vertex  $(w_1,\cdots ,w_n)$ in  $\overline{Q_{B}}$,
 \begin{itemize}
 \item[(a)] we denote  $d_{n}^{0}$ the arrow
 $$\xymatrix{ (w_{1},w_2, \cdots ,  w_{n}) \ar[rd]^{d_{n}^{0}}_{w_1\ot 1} &\\
   & (w_{2},\cdots, ,w_{n})}$$
   with weight $w_1\ot 1$, which is put under the arrow (we will do the same in the following);
 \item[(b)] for $1\leq i\leq n-1$,  assume that $w_iw_{i+1}\equiv \sum_j \lambda_j u_j\  \mathrm{mod}\ I$ with all $u_j\in \NonTip(I)\backslash Q_0$ and $\lambda_j\in k^*=k\setminus \{0\}$;  by abuse of notations, we still denote by $d_{n}^i$  all arrows of the form:
 $$\xymatrix{ (w_{1},\cdots , w_{i-1}, w_{i}, w_{i+1}, w_{i+2}, \cdots, w_{n}) \ar[rd]^{d_{n}^{i}}_{(-1)^i \cdot \lambda_j \otimes 1} &\\
   & (w_{1},\cdots,w_{i-1}, u_j ,w_{i+2}, \cdots,w_{n});}$$
 \item[(c)]  the arrow
$$\xymatrix{ (w_{1}, \cdots, w_{n-1}, w_{n}) \ar[rd]^{d_{n}^{n}}_{(-1)^n 1\ot w_n} &\\
   & (w_{1},\cdots, ,w_{n-1})}$$
   with weight $(-1)^n 1\ot w_n$ is denoted by $d_n^{n}$. Note that the weight here should be written as $(-1)^n 1\ot w_n^{op}$, but for simplicity, we just write it as $(-1)^n 1\ot w_n$.
 \end{itemize}

\begin{Rem} \label{Def:new-weighted-quiver-overline}
We can define all the concepts in Section~\ref{Sect: morse theory} for $\overline{Q_{B}}$ analogously. Notice that the only difference between $Q_{B}$ and $\overline{Q_{B}}$ is  that there exist parallel (or multiple)  arrows in $\overline{Q_{B}}$ whose weights sum to that of the arrow with the same starting vertex and target in $Q_B.$ It follows that  Theorem~\ref{Thm: main result of algebraic Morse theory} yields the same result if we use $\overline{Q_{B}}$ instead of $Q_{B}.$  Note that the notion $\overline{Q_{B}}$ and the notion $\overline{Q_B}^\M$ below in the present paper is only valid for the reduced two-sided bar resolution $B(A,A)$.
\end{Rem}

For $w\in\mathcal{B}$, let $V_{w,i}$ be the vertices $(w_1,\cdots,w_n)$ in $\overline{Q_B}$ such that $w=w_1\cdots w_n$ and $i$ is the largest integer $i\geq -1$ such that $(w_1,\cdots,w_{i+1})$ is an $i$-chain. Let $V_w=\bigcup_iV_{w,i}$. Thus $(w_1,\cdots,w_n)\in V_{w,-1}$ if and only if $w_1\notin Q_1$, $(w_1,\cdots,w_n)\in V_{w,0}$ if and only if $w_1\in Q_1$ and $(w_1,w_2)$ is not a $1$-chain, etc.

We define a partial matching $\M$ to be the set of arrows of the following form in $\overline{Q_B}$:
\begin{eqnarray}\label{arrow} (w_1,\cdots,w_{i+1},w_{i+2}',w_{i+2}'',w_{i+3},\cdots,w_n) \to (w_1,\cdots,w_{i+1},w_{i+2},w_{i+3},\cdots,w_n),\end{eqnarray}
where $(w_1,\cdots,w_n)\in V_{w,i}$, $w_{i+2}' w_{i+2}''=w_{i+2}$ and $(w_1,\cdots,w_{i+1},w_{i+2}',w_{i+2}'',w_{i+3},\cdots,w_n)\in V_{w,i+1}$.
 Note that in this case $i$ is necessarily less than $n-1$.  Indeed, $\M$ is a partial matching: clearly no vertex in $\M$ is the starting vertex of more than one edge in $\M$; if a vertex $(w_1,\cdots,w_{i+1},w_{i+2},w_{i+3},\cdots,w_n)$ in $\M$ were the target vertex of more than one edge in $\M$, then there would have to be two different decompositions $w_{i+2}=w_{i+2}' w_{i+2}''=v_{i+2}' v_{i+2}''$ such that $w_{i+1}w_{i+2}'$ has a unique factor in $W$ which is a  right factor  of $w_{i+1}w_{i+2}'$ and also that $w_{i+1}v_{i+2}'$ has a unique factor in $W$ which is a  right factor  of $w_{i+1}v_{i+2}'$, and this is a contradiction; the situation $$(w_1,\cdots,w_{i+1},w_{i+2}',w_{i+2}'',w_{i+3},\cdots,w_n) \to (w_1,\cdots,w_{i+1},w_{i+2},w_{i+3},\cdots,w_n)\in V_{w,i},$$
$$(w_1,\cdots,w_{i+1},w_{i+2}',w_{i+2}^{(3)},w_{i+2}^{(4)},w_{i+3},\cdots,w_n) \to (w_1,\cdots,w_{i+1},w_{i+2}',w_{i+2}'',w_{i+3},\cdots,w_n)\in V_{w,i+1}$$
cannot occur since this would imply $w_{i+2}$ lies in $\langle \Tip(I)\rangle$; moreover, the arrow  (\ref{arrow})  in $\M$  represents an invertible homomorphism since its weight is $(-1)^{i+2}$.

Recall from the last section that the above partial matching gives a new decorated quiver $\overline{Q_B}^\M$ from the original decorated quiver $\overline{Q_B}$ by reversing the arrows in $\M$.
  The arrows $d_n^i$ in $\overline{Q_B}$ which remain unchanged in $\overline{Q_B}^\M$ will be drawn by thick arrows, still denoted by $d_n^i$;  an arrow $d_n^i$ lying in $\M$ will be drawn in dotted arrows with reverse direction, denoted by $d_{n}^{-i}$.
   By definition, the arrow $d_n^i$  in $\M$ necessarily has   $1\leq i\leq n-1$ and $w_iw_{i+1}\in \NonTip(I)\backslash Q_0$, so the arrow $d_{n}^{-i}$ has the form
   $$\xymatrix{   & (w_{1},\cdots,w_{i-1}, w_iw_{i+1} ,w_{i+2}, \cdots,w_{n})\ar@{..>}[ld]_{d_{n}^{-i}}^{(-1)^{i+1} }  \\
  (w_{1},\cdots , w_{i-1}, w_{i}, w_{i+1}, w_{i+2}, \cdots, w_{n}) & .}$$


The following result should be compared with   \cite[Lemma 9]{S06}.
\begin{Thm}  \label{thm:skoeldberg}
The partial matching $\M$ defined above is a Morse matching of $\overline{Q_B}$ such that the set of critical vertices in $n$-th component is identified with the set $W^{(n-1)}$ of $n-1$-chains.
\end{Thm}

\begin{Proof}
First we identify the critical vertices in the $n$-th component. Suppose $(w_1,\cdots,w_n)\in V_{w,i}$ is a critical vertex. Then $(w_1,\cdots,w_{i+1})$ is an $i$-chain. Suppose $i<n-1$. There are two cases: $w_{i+1} w_{i+2}\in \langle \Tip(I)\rangle$ or $w_{i+1} w_{i+2}\notin \langle \Tip(I)\rangle$. When $w_{i+1} w_{i+2}\in \langle \Tip(I)\rangle$, there is a decomposition $w_{i+2}=w_{i+2}' w_{i+2}''$ with $w_{i+2}'$ minimal such that $w_{i+1}w_{i+2}'\in \langle \Tip(I)\rangle$. Since $(w_1,\cdots,w_{i+2})$ is not an $i+1$-chain, $w_{i+2}''$ has nonzero length, which means that there is an edge $(w_1,\cdots,w_{i+1},w_{i+2}',w_{i+2}'',\cdots,w_n) \to (w_1,\cdots,w_n)$ in $\M$. This contradicts that $(w_1,\cdots,w_n)$ is a critical vertex. When $w_{i+1} w_{i+2}\notin \langle \Tip(I)\rangle$, $(w_1,\cdots,w_{i+1}w_{i+2},\cdots,w_n)\in V_{w,i-1}$ and there is an edge $(w_1,\cdots,w_n) \to (w_1,\cdots,w_{i+1}w_{i+2},\cdots,w_n)$ in $\M$. This is also a contradiction. Hence $i=n-1$ and $(w_1,\cdots,w_{n})\in  W^{(n-1)}  $. On the other hand, it is obvious that all vertices $(w_1,\cdots,w_n)$ in  $W^{(n-1)}$   are critical vertices in $n$-th component. Thus the set of critical vertices in the $n$-th component is identified with the set $W^{(n-1)}$ of $n-1$-chains.

 Next, we show that any zigzag path in $\overline{Q_{B}}^{\M} $ is of finite length. Consider a vertex $(w_1,\cdots,w_n) \in V_w$, and look at the corresponding differential
$$d([w_1 | \cdots |w_n])=w_1 [w_2 |\cdots| w_n]+\sum\limits_{i=1}^{n-1} (-1)^i [w_1| \cdots|w_i w_{i+1}|\cdots |w_n]+(-1)^n [w_1| \cdots | w_{n-1}]w_n.$$
The element $w_1 [w_2 |\cdots| w_n]$ (resp. $[w_1| \cdots | w_{n-1}]w_n$) is in the component corresponding to the vertex $(w_2,\cdots,w_n)$ (resp. $(w_1,\cdots,w_{n-1})$), and $w_2\cdots w_n\prec w_1\cdots w_n$ (resp. $w_1\cdots w_{n-1}\prec w_1\cdots w_n$). The elements $[w_1| \cdots|w_i w_{i+1}|\cdots |w_n]$ can all be written as linear combinations of elements in components corresponding to $(w_1,\cdots,w_{i-1},u,w_{i+2},\cdots,w_n)$, where $w_1\cdots w_{i-1}uw_{i+2}\cdots w_n\preceq w_1\cdots w_n$, with equality or inequality depending on whether $w_{i} w_{i+1}\notin \langle \Tip(I)\rangle$ (this is the only case such that the vertex $(w_1,\cdots,w_{i-1},u,w_{i+2},\cdots,w_n)$ remains in $V_w$ with $u=w_{i} w_{i+1}$) or not (in this case $w_{i} w_{i+1}\in \Tip(I)$ and $w_{i} w_{i+1}=\sum_j\lambda_ju_j$ (mod $I$) with all $u_j\prec w_{i} w_{i+1}$ or $w_{i} w_{i+1}=0$).
So for a thick arrow $v\to v'$ in $\overline{Q_B}^\M$ with $v\in V_{w} $ and $v'\in V_{w'},$ we have $w\succeq w'.$   On the other hand, for a dotted arrow $v\dashrightarrow v'$ in $\overline{Q_B}^\M$,  we have $v, v'\in V_{w}$ for some $w\in\mathcal{B}$ by the definition of $\M$.

 Let $p=\alpha_1   \alpha_2  \cdots $ be a zigzag path in $\overline{Q_B}^\M$. We are going to prove that $p$ has finite length. By the well-ordering of $\preceq$ and the observations above, without loss of generality we may assume that the starting vertex of each arrow in $p$ belongs to $V_w$ with the same $w\in \mathcal{B}.$ Let
$$v_k \overset{\alpha_k}{\dashrightarrow} v_{k+1} \overset{\alpha_{k+1}}{\to} v_{k+2}$$ be the segment of $p$ of length $2.$ By the construction of $\M,$ the dotted arrow $\alpha_k$ has the form
 $$d_{n+1}^{-(i+2)}:v_k=(w_1,\cdots,w_{i+1},w_{i+2},w_{i+3},\cdots,w_n) \dashrightarrow v_{k+1}=(w_1,\cdots,w_{i+1},w_{i+2}',w_{i+2}'' ,w_{i+3},\cdots,w_n)$$
 with $v_k \in V_{w,i}$ and $v_{k+1}\in V_{w,i+1},$ $-1\leq i \leq n-2.$ As $v_{k+2}\in V_{w}, $ the thick arrow $\alpha_{k+1}\neq d_{n+1}^{l},0\leq l\leq i+1.$ Notice that $\alpha_k=d_{n+1}^{-(i+2)},$ thus the arrow $d_{n+1}^{i+2}$ with starting vertex $v_{k+1}$ does not exist in $\overline{Q_B}^\M.$ Hence $\alpha_{k+1}=d_{n+1}^{l}$ with $i+3\leq l \leq n+1$ which shows that the first $i+2$ components of $v_{k+1}$ and $v_{k+2}$ coincide. Then we have $v_{k+2} \in V_{w,i+j}$ with $j>0$. As the length of  $w\in \mathcal{B}$ is finite, the subscript $i$ of $V_{w,i}$ has a finite upper bound. So the zigzag path $p$ is of finite length.
 Hence $\M$ is a Morse matching in terms of Proposition~\ref{Prop: sufficient condition for Morse matching}.

\end{Proof}

Hence by Theorem~\ref{Thm: main result of algebraic Morse theory},  the reduced two-sided bar resolution $B(A,A)$ of the algebra $A$ is homotopy equivalent to a complex $(B(A,A)^\M,d^{\M})$ associated to the quiver $\overline{Q_B}^\M$.
We will call $ (B(A,A)^\M,d^{\M})$ a two-sided Anick resolution of $A$ according to \cite[Theorem 4]{S06}, and this resolution can be described as follows:
  for $n\geq 0$, the $n$-th component is $A\otimes_E kW^{(n-1)}\otimes_EA$, and the differential from the $n$-th to the $(n-1)$-th component corresponds to the sum
         $$ \sum_{\substack{ p\in \Path^\M(w, w') \\ w\in W^{(n-1)}, w'\in W^{(n-2)}} } \varphi_p^\M$$
    of all zigzag paths in (LFH1). Theorem~\ref{Thm: main result of algebraic Morse theory} has the following obvious corollary.

\begin{Cor} If $W$ is a finite set and the Ufnarovski\u\i~graph $Q_W$
has no oriented cycles, that is, $Q_W$ is a finite acyclic quiver, then the algebra $A$ has finite global dimension.
\end{Cor}

\begin{exam} \label{exam:two-sided-Anick-resolution} A concrete calculation shows that the rightmost part of the two-sided Anick resolution of the algebra $A$ in Example \ref{exam:Uf-graph} is the following:
\begin{center}
	\begin{tikzcd}
	\cdots \arrow{r} & P_2 \arrow{r}{d^{\M}_2} & P_1 \arrow{r}{d^{\M}_1} & P_0 \arrow{r}{\epsilon} & A \arrow{r} & 0,
	\end{tikzcd}
\end{center}
 where $$P_n:=A\otimes_E kW^{(n-1)}\otimes_EA\cong \bigoplus_{(w_1,\cdots,w_n)\in W^{(n-1)}} A^e \cdot [w_1 | \cdots | w_n],$$
$$W^{(-1)}=\{e_1,e_2,e_3\}, W^{(0)}=\{(a),(b),(a'),(b')\}, W^{(1)}=\{(a,b),(b,b'),(b',a'),(a,a'a),(a',aa')\},$$
$\epsilon$ is the multiplication,  $d^{\M}_1$ is given by the matrix (here the element $b\otimes 1$ means that the generator $(b)$ in $P_1$ maps to $b\otimes 1$ times the generator $e_3$ in $P_0$, etc.)
$$\left(\begin{array}{cccc} -1\otimes a & a\otimes 1&0 \\ 0& -1\otimes b& b\otimes 1 \\ a'\otimes 1 & -1\otimes a'& 0 \\ 0& b'\otimes 1 & -1\otimes b' \end{array} \right),$$
$d^{\M}_2$ is given by the matrix
$$\left(\begin{array}{ccccc} 1\otimes b & a\otimes 1&0 & 0 \\ -a'\otimes 1& 1\otimes b'&-1\otimes a&b\otimes 1 \\ 0&0&b'\otimes 1 & 1\otimes a' \\ 1\otimes a'a+aa'\otimes 1& 0&a\otimes a &0\\ a'\otimes a'& 0& 1\otimes aa'+a'a\otimes 1&0 \end{array} \right).$$
For example, all the zigzag paths from the vertex $(b,b')$ can be calculated using the following diagram:
$$
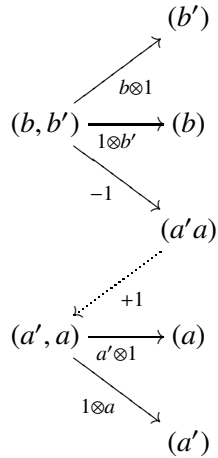
\xymatrix{& (b') \\
 (b,b')\ar[ru]_{b\otimes 1}\ar[r]_{1\otimes b'}\ar[rd]_{-1} & (b)\\
 & (a'a)\ar@{..>}[ld]^{+1}\\
 (a',a)\ar[r]_{a'\otimes 1}\ar[rd]_{1\otimes a}& (a)\\
 & (a') }$$
\end{exam}

\bigskip

\underline{}\section{Minimality criterion for  two-sided Anick resolutions}\label{Sect: minimality}

In this section we will consider in which case a two-sided Anick resolution is minimal.

\begin{Def} \label{Def: Minimal resolution definition}
Let $A=kQ/I$ be an algebra as in Section~\ref{Sect: Two-sided Anick resolutions}, $E=\oplus_{e\in Q_0}ke$ be its semisimple subalgebra. Then $A$ is an augmented $E$-algebra with augmented ideal $A_+$.

A projective resolution $(P_*,d_*)$ of  a left $A$-module $M$ is minimal if the induced map $$1 \otimes d_*: E \otimes_{A} P_* \to E \otimes_{A} P_{*-1}$$ is zero. A projective resolution $(P'_*,d'_*)$ of  an $A$-bimodule $M'$ is minimal if the induced map $$1 \otimes d'_*\otimes 1: E \otimes_{A} P'_* \otimes_{A} E \to E \otimes_{A} P'_{*-1}\otimes_{A} E$$ is zero.
\end{Def}

 Our definition of minimality is consistent with the usual one in the literature, see, for example \cite[Page 325]{AG}.

We keep the notations  in Section~\ref{Sect: Two-sided Anick resolutions}.

\begin{Def}\label{Def: Reduction}
Let $w, w'\in \mathcal{B}$.  Define a reduction step from $w$ to $w'$ with coefficient $\lambda\in k^{*}$, denoted by $w \Longrightarrow_{\lambda} w'$,  if there exist $u, v \in \mathcal{B}$ and $f\in I$ such that
 \begin{itemize}
\item[(a)] $\Tip(f) \in W$; 
\item[(b)] $w = \Tip(ufv)=u \Tip(f) v;$
\item[(c)] $-\lambda w'= u p v$ where  $p\neq \Tip(f)$ is a monomial appearing in $f$.
 \end{itemize}
 We say $w$ converges to $w'$ if there is a sequence of reduction steps $w \Longrightarrow_{\lambda_1}  u_1  \Longrightarrow_{\lambda_2}   \cdots  \Longrightarrow_{\lambda_m}u_m \Longrightarrow_{\lambda_{m+1}} w'$ with $u_1, \cdots , u_m \in \mathcal{B}$.
\end{Def}







 \begin{Rem}\label{Remark: arrows vs rewriting}
 For  a thick arrow
$$ (u_1, \cdots, u_n) \overset{d_{n}^i}{\to}  (v_1, \cdots, v_{n-1})$$
 with $1\leq i \leq n-1$, we have either $u=u_1\cdots u_n \Longrightarrow_{\lambda} v=v_1\cdots v_{n-1}$ or $u=v$. Similarly,
for a dotted arrow
$$\xymatrix{    & (u_{1},\cdots , u_{n-1}) \ar@{..>}[ld]_{d_{n}^{-i}}\\
  (u_{1}, \cdots, u_{i-1}, u_{i}', u_{i}'', u_{i+1}, \cdots, u_{n-1})&}$$
with $1\leq i \leq n-1$, we have     $u_1\cdots u_{n-1} = u_1\cdots u_{i-1}u_i'u_i'' u_{i+1}\cdots u_n$.
 \end{Rem}

The following result says that if an arbitrary $(n-1)$-chain cannot converge to an $(n-2)$-chain, then the two-sided Anick resolution is minimal, thus providing a handy minimality criterion.

\begin{Thm}[Minimality criterion]\label{Thm: minimal criterion}
 With the notations in Section~\ref{Sect: Two-sided Anick resolutions}. For an arbitrary $(n-1)$-chain $(w_1, \cdots , w_{n})$ with $n\geq 1$, if $w=w_1\cdots w_n$ cannot converge to $u=u_1\cdots u_{n-1}$ for any $(u_1,\cdots, u_{n-1}) \in W^{(n-2)}$,  then  the two-sided Anick resolution $ (B(A,A)^\M,d^{\M})$ is minimal.
 \end{Thm}

\begin{Proof}

It is clear that the two-sided Anick resolution $ (B(A,A)^\M,d^{\M})$ is minimal if and only if
$$\mrm{Im} (d_{n}^{\M}) \subseteq (A^{e})_{+} \otimes_{E^{e}} B(A,A)_{n}^\M$$
 for $n\geq 0.$
So it suffices to prove that each zigzag path from $ (w_1,\cdots, w_n)\in W^{(n-1)}$ to $ (w'_1,\cdots ,w'_{n-1}) \in W^{(n-2)}$ in $\overline{Q_B}^\M$ has weight $\lambda \in (A^{e})_{+}$ for $n\geq 1.$
Let $p=\alpha_1\alpha_2\cdots \alpha_m$ be such a zigzag path.
By Remark~\ref{Remark: arrows vs rewriting},  if the thick arrows appearing in $p$   have the form $d_{n}^{i},$ $1\leq i\leq n-1,$ necessarily  $w=w_1\cdots w_n$ converges to $w'=w_1'\cdots w_{n-1}'$. Thus, by our assumption, there exists a thick arrow $\alpha_j=d_{n}^{0}$ or $\alpha_j=d_{n}^{n} $ whose weight lies in $(A^e)_+.$
As the weight of $p$ is the product of the weight of all arrows appearing in $p.$
Hence the weight of $p$  lies in $(A^e)_+$.

\end{Proof}

 In practice, when the Gr\"{o}bner-Shirshov basis is homogeneous, it is sometimes helpful to use Theorem~\ref{Thm: minimal criterion} to determine whether the two-sided Anick resolution is minimal  by calculating the degrees of the Anick chains; for  example, if each $n$-chain has different degree with any $(n-1)$-chain, it immediately follows that the two-sided Anick resolution is minimal.

\begin{Rem}\label{Remark: minimal at 012} Notice that the hypothesis of the above result  holds unconditionally  for  $n= 0, 1, 2,$  so the two-sided Anick resolution is minimal in degree $0, 1, 2.$

\end{Rem}

\begin{Rem} The original papers of Anick \cite{Anick}  and Anick-Green \cite{AG} only consider one-sided Anick resolutions.
    It is easy to see that  the minimality criterion Theorem~\ref{Thm: minimal criterion} still works for one-sided Anick resolutions.
    \end{Rem}

\begin{Rem} In the article \cite{JW}, the authors gave a sufficient condition for the minimality of   one-sided Anick resolutions:
\begin{itemize}
\item  There does not exist an arrow of type II  in $\overline{Q_B}^\M$, that is,  an arrow of the form
$$d_n^i: (w_1, \cdots , w_n) \to  (w_1, \cdots ,w_{i-1}, v ,w_{i+2}, \cdots ,w_n)$$
with $w_iw_{i+1} \Longrightarrow_{\lambda} v $ and $(w_1, \cdots ,w_{i-1}, v ,w_{i+2}, \cdots ,w_n)\in W^{(n-2)}$ (Note that they do not ask
$(w_1, \cdots , w_n)$ to be an $(n-1)$-chain).
\end{itemize}

However,  there is a counter-example.

Let  $A=k\langle x_1,x_2,\cdots , x_7\rangle /I$  with $I=\langle x_1x_2x_3-x_6x_7, x_3x_4x_5, x_6x_7x_4x_5\rangle$.  Fix the order  $x_1>x_2>\cdots >x_7$.  The set of the reduced monomial generators of $\langle\Tip(I)\rangle$ with respect to the left length-lexicographic order  is $W=\{ x_1x_2x_3, x_3x_4x_5 ,x_6x_7x_4x_5 \}$. We list all the $i$-chains ($i\geq -1$) as follows:
\begin{itemize}
\item[(a)] $W^{(-1)}$ is a singleton;
\item[(b)] $W^{(0)} = \{ x_1, \cdots, x_7 \}$;
\item[(c)] $W^{(1)}=\{ (x_1, x_2x_3), (x_3,x_4x_5), (x_6, x_7x_4x_5) \}$;
\item[(d)] $W^{(2)}= \{ (x_1, x_2x_3, x_4x_5) \}$;
\item[(e)] $W^{(i)}=\emptyset, i \geq 3.$
\end{itemize}
 The (one-sided) Anick resolution of $k$ has the form $$0 \to A\otimes W^{(2)} \overset{d^{\M}_3}{\to}   A\otimes W^{(1)}  \overset{d^{\M}_2}{\to}  A\otimes W^{(0)} \overset{d^{\M}_1}{\to} A \to k \to 0.$$
We now compute the differential $d_3^{\M}$.  There are exactly two zigzag paths starting from a $2$-chain to a $1$-chain as follows:
$$(a) 
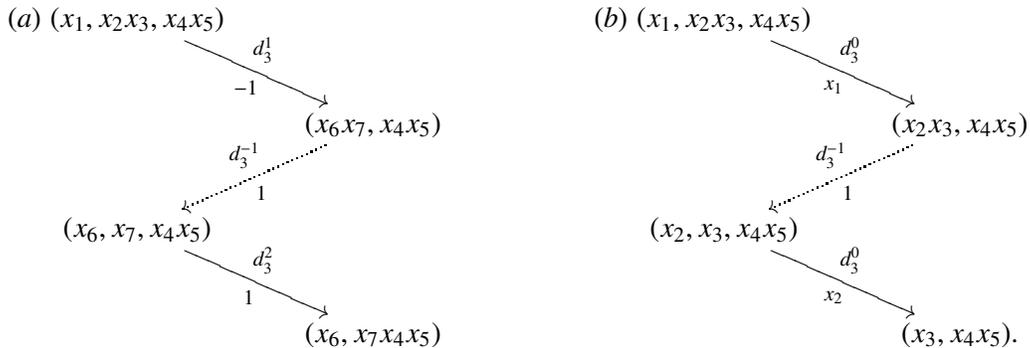
\xymatrix{ (x_1, x_2x_3, x_4x_5)\ar[rd]^{d_{3}^{1}}_{-1}&\\
  &(x_6x_7, x_4x_5)\ar@{..>}[ld]_{d_3^{-1}}^{1}  \\
   (x_6, x_7,x_4x_5) \ar[rd]^{d_3^2}_{1}&\\
 &  (x_6, x_7x_4x_5)   }
\hspace{5em} (b) \xymatrix{ (x_1, x_2x_3, x_4x_5)\ar[rd]^{d_{3}^{0}}_{x_1}&\\
  &(x_2x_3, x_4x_5)\ar@{..>}[ld]_{d_3^{-1}}^{1}  \\
   (x_2, x_3, x_4x_5) \ar[rd]^{d_3^0}_{x_2}&\\
 &  (x_3, x_4x_5)   .}$$
 We have  $d_3^{\M}(1\otimes (x_1, x_2x_3, x_4x_5))= -1\otimes (x_6, x_7x_4x_5) + x_1x_2\otimes (x_3,x_4x_5)$. It follows that the induced map $1\otimes d_3^{\M}\neq 0$.   So the resolution is not minimal.

However,  since  all reduction steps are induced by  the reduction $x_1x_2x_3 \Longrightarrow_1 x_6x_7$,  no reduction of type $II$ appears
because  there is no    Anick chain of the form $(w_1, \cdots, ux_6x_7v, \cdots, w_n)$.





\end{Rem}

\begin{Rem}
In general, the criterion in Theorem~\ref{Thm: minimal criterion} is NOT a necessary condition.

Let $B=k\langle x_1, x_2, x_3, x_4, x_5 \rangle  / I$ with $I=\langle x_1x_2x_3- x_1x_5, x_2x_3x_4-x_5x_4, x_1x_5x_4\rangle $. Fix the order $x_1>x_2>x_3>x_4>x_5$.
The set of the reduced monomial generators of $\langle\Tip(I)\rangle$ with respect to the left length-lexicographic order  is $W=\{ x_1x_2x_3, x_2x_3x_4 , x_1x_5x_4 \}$. We list all the $i$-chains ($i\geq -1$) as follows:
\begin{itemize}
\item[(a)] $W^{(-1)}$ is a single point set;
\item[(b)] $W^{(0)} = \{ x_1, \cdots, x_5 \}$;
\item[(c)] $W^{(1)}=\{ (x_1, x_2x_3), (x_2,x_3x_4), (x_1, x_5x_4) \}$;
\item[(d)] $W^{(2)}= \{ (x_1, x_2x_3, x_4) \}$;
\item[(e)] $W^{(i)}=\emptyset, i \geq 3.$
\end{itemize}
 The  two-sided  Anick resolution of $A$ has the form
 $$0 \to A\otimes W^{(2)}  \otimes A \overset{d^{\M}_3}{\to}   A\otimes W^{(1)}\otimes A  \overset{d^{\M}_2}{\to}  A\otimes W^{(0)}\otimes A \overset{d^{\M}_1}{\to} A\otimes A \to A \to 0.$$
 For the $2$-chain $(x_1, x_2x_3, x_4)$, there are two reduction steps converging to the $1$-chain $(x_1, x_5x_4)$:  $$w=\underline{x_1x_2x_3}x_4\Longrightarrow_1 x_1x_5x_4 \text{ and } w=x_1\underline{x_2x_3x_4}\Longrightarrow_{1} x_1x_5x_4.$$
 So the assumption of Theorem~\ref{Thm: minimal criterion} is not fulfilled.

 However,  this   two-sided  Anick resolution of $A$ is minimal.
 For the differential $d_3^{\M}$, there are  four  zigzag paths starting from a $2$-chain to an $1$-chain as follows:
  $$ (a) 
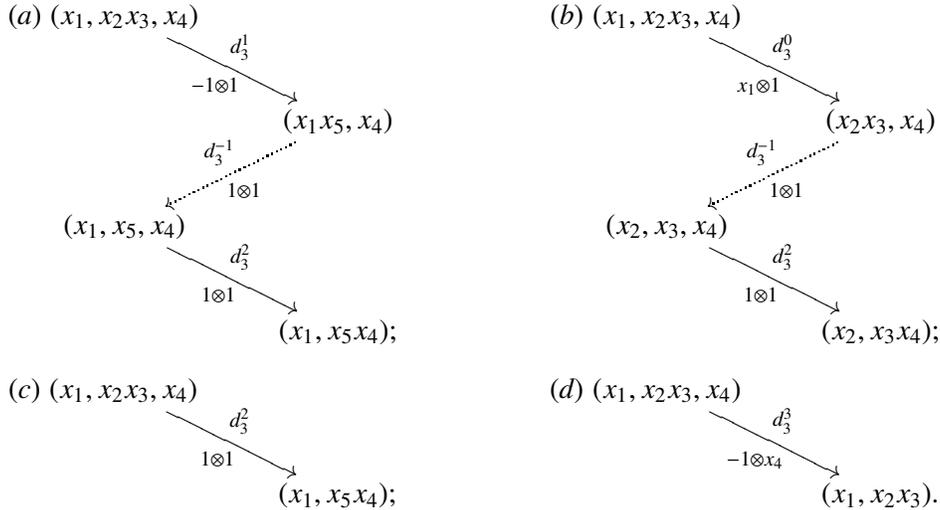
\xymatrix{ (x_1, x_2x_3, x_4)\ar[rd]^{d_{3}^{1}}_{-1\otimes 1}&\\
  &(x_1x_5, x_4)\ar@{..>}[ld]_{d_3^{-1}}^{1\otimes 1}  \\
   (x_1, x_5,x_4) \ar[rd]^{d_3^2}_{1\otimes 1}&\\
 &  (x_1, x_5x_4);   }
 \hspace{5em}
 (b) \xymatrix{ (x_1, x_2x_3, x_4)\ar[rd]^{d_{3}^{0}}_{x_1\otimes 1}&\\
  &(x_2x_3, x_4)\ar@{..>}[ld]_{d_3^{-1}}^{1\otimes 1}  \\
   (x_2, x_3, x_4) \ar[rd]^{d_3^2}_{1\otimes 1}&\\
 &  (x_2, x_3x_4)   ;}
 $$
 $$(c)\xymatrix{ (x_1, x_2x_3, x_4)\ar[rd]^{d_{3}^{2}}_{1\otimes 1}&\\
  &(x_1, x_5x_4)  ;} \hspace{5em}
 (d)\xymatrix{ (x_1, x_2x_3, x_4)\ar[rd]^{d_{3}^{3}}_{-1\otimes x_4}&\\
  &(x_1, x_2x_3)  .}$$
 We thus have
 \begin{eqnarray*}
 &&d_3^{\M}(1\ot (x_1,x_2x_3,x_4)\otimes 1)\\
 &&= -1\ot(x_1,x_5x_4)\otimes 1 + 1\ot (x_1,x_5x_4) \otimes 1 +x_1\ot (x_2,x_3x_4)\otimes 1- 1\otimes (x_1,x_2x_3)\otimes x_4\\
 &&= x_1\ot (x_2,x_3x_4)- 1\otimes (x_1,x_2x_3)\otimes x_4.
 \end{eqnarray*} So the induced map $1_k\ot d_3^{\M}\otimes 1_k=0$.   By Remark~\ref{Remark: minimal at 012},  $1_k\ot d_2^{\M}\otimes 1_k$ and $1_k\ot d_1^{\M} \otimes 1_k$ vanish  as well. Hence the  two-sided  Anick resolution is minimal.
\end{Rem}

\bigskip

\section{Homological properties of Chinese algebras}

 In this section, we will use the two-sided Anick resolution to study  the homological properties of   Chinese algebras of finite rank \cite{CEKNH, CQ}.

Let $X= \{ x_1, \cdots, x_n\}$ with $n\geq 1,$
 $X^*$ be the free monoid generated by $X$. The Chinese congruence is the congruence on  $X^*$ generated by  the following relations:
\begin{itemize}
\item[(a)] $x_i x_j x_k=x_i x_k x_j= x_j x_i x_k,  \forall i >j > k$,
\item[(b)] $x_i x_j x_j =x_j x_i x_j$, $x_i x_i x_j= x_i x_j x_i, \forall i > j$.
\end{itemize}
The Chinese monoid $\CH(X)$ (of rank $n$) is the quotient monoid of the free monoid $X^*$ by the Chinese congruence.

\begin{Def}[Chinese algebra]
Let $k$ be a field. The Chinese algebra $A$ (of rank $n$) is the semigroup algebra of the Chinese monoid $\CH(X)$.
\end{Def}

Equivalently, $A$ is the algebra with relation by $k\langle X|T\rangle =k\langle x_1,\cdots ,x_n\rangle / I$ with
$$I=\langle x_i x_j x_k-x_j x_i x_k, x_i x_k x_j- x_j x_i x_k, x_i x_j x_j -x_j x_i x_j, x_i x_i x_j- x_i x_j x_i \rangle_{i>j>k}.$$

Y. Chen and J. Qiu obtained the Gr\"{o}bner-Shirshov basis $\mathcal{G}$ with respect to the left length-lexicographic order on $X^*$ generated by $$x_n > x_{n-1} > \cdots >x_1$$
as follows:
\begin{Thm}[\cite{CQ}]
    The Gr\"{o}bner-Shirshov basis $\mathcal{G}$ for the ideal $I$ with respect to the left length lexicographic order on $X^*$ consists of
\begin{itemize}
\item[(a)] $x_ix_jx_k-x_jx_ix_k, x_ix_kx_j-x_jx_ix_k$, $\forall \ i>j>k$;
\item[(b)] $x_ix_jx_j-x_jx_ix_j,x_ix_ix_j-x_ix_jx_i$, $\forall \ i>j$;
\item[(c)] $x_ix_jx_ix_k-x_ix_kx_ix_j$, $\forall \ i>j>k$.
\end{itemize}

\end{Thm}

\noindent  Note that the above Gr\"{o}bner-Shirshov basis of the Chinese algebra $A$ is homogeneous and the set of reduced monomial generators of $\Tip(I)$ is
$$W=\Tip(\mathcal{G})=\{ x_ix_jx_k, x_ix_kx_j, x_ix_jx_ix_k\}_{i>j>k }\cup \{ x_ix_jx_j, x_ix_ix_j \}_{i>j}.$$
By Theorems~\ref{thm:skoeldberg} and \ref{Thm: main result of algebraic Morse theory}, we obtain a free bimodule resolution $(P_*,  d_{*} )$ of $A$  (for simplicity here and in the following we often write $d_{*}^{\M}$ as $d_{*}$)  with  $$P_i=A\otimes_k kW^{(i-1)} \otimes_k A, i\geq 0.$$

The following lemma implies that  the length of the resolution $(P_*,  d_{*} )$ is equal to $\frac{n(n+1)}{2}$.

 \medskip

 \begin{Lem}\label{Lem: chain structure of Chinese alg}
 With the notation above,  the sets of Anick chains $W^{(i)}$ have the  following properties:
\begin{itemize}
\item[(a)] $|W^{(\frac{n(n+1)}{2}-1)}|=1$ and  the unique element   in $W^{(\frac{n(n+1)}{2}-1)}$  is of  maximal degree  among all the  Anick chains;
\item[(b)] $W^{(k)}=\emptyset$, $\forall k\geq \frac{n(n+1)}{2}$;

\item[(c)] the Ufnarovski\u\i~graph  $Q_W$ (cf. Section 4) is a finite quiver  and it  has no oriented cycle;

\item[(d)] each $W^{(i)}$ is a finite set.
\end{itemize}

\end{Lem}

\begin{Proof}
We prove (a)(b)  by induction on   $n$,  the rank of the Chinese algebra $A$.

The case $n=1$ is trivial as   $W=\emptyset$ and there are only two Anick chains: one $(-1)$-chain and one $0$-chain $x_1$.

For the case $n=2$,   we have $W=\{x_2 x_2 x_1, x_2x_1x_1  \}$. It follows  that $W^{(2)}=\{(x_2,x_2x_1,x_1)\},$ $W^{(1)}=\{(x_2,x_2x_1),(x_2,x_1x_1)  \},$ $W^{(0)}=\{x_1,x_2\},$ $W^{(-1)}=\{ \ast \}$ and $W^{(k)}=\emptyset$ for $k\geq 3.$
Thus (a)(b) hold for $n=2$.

 Inductively,  let $w=(w_1, \cdots, w_m)$ be  an  Anick chain of the maximal weight.
It is clear that $w_1\in \{ x_1,\cdots, x_n\}$. If $w_1\neq x_n$, we can construct a new Anick chain $(x_n, x_n x_{n-1}, w_1, \cdots, w_m)$  whose weight and degree are greater   than those of $w$.  So we have $w_1=x_n$.

 It follows from the definition of Anick chains that $w_1 w_2 =x_nw_2\in W$,  then we have
$$w_2 \in \{x_j x_k, x_kx_j, x_j x_n x_k  \}_{k<j<n}  \cup \{ x_j x_j, x_n x_j \}_{j<n}. $$
Let's  discuss these cases  except for the case $w_2=x_nx_{n-1}.$
\begin{itemize}
 \item[(i)]    $w_2=x_j x_k$   with $k<j<n$.
 One can see that  the sequence $(x_n, x_n x_{n-1}, x_j, w_2, w_3,\cdots, w_m)$ is an Anick chain whose   weight and degree    are greater than those of $w$.

 \item[(ii)]   $w_2=x_kx_j$  with $k<j<n$. As there is no monomial beginning with $x_k x_j$  in $W$,   we have $x_jw_3\in W$. Thus, the weight and degree of the Anick chain $(x_n, x_n x_{n-1}, x_j, w_3,\cdots,w_m)$ is greater than those of $w$.

 \item[(iii)]    $w_2= x_j x_n x_k$ with $k<j<n.$   Notice that $x_j x_n x_k w_3 \notin W$ for any $w_3\in X^{*},$ so we have  either $x_nx_kw_3\in W$  or $x_kw_3\in W.$
     If $x_nx_kw_3\in W$,  the weight  and degree of the Anick chain $(x_n,x_nx_{n-1}, x_nx_{n-2}, w_3,\cdots ,w_m)$ are greater  than those of $w$.
    Now we assume that $ x_kw_3\in W.$  It follows that the degree of $w_3$ is $2$ or $3$ and $w_3=x_lw_3'$ with $l\leq k.$ Thus we have $x_nx_kx_l\in W$ which contradicts the fact that $w$ is an Anick chain.

 \item[(iv)]  $w_2=x_j x_j$  with $j<n.$ We have either  $x_jx_jw_3\in W$ or  $x_jw_3\in W$. If $x_jx_jw_3\in W$, then $w_3=x_{\ell}$ with $\ell<j$. Thus the weight and degree of the Anick chain $(x_n,x_nx_{n-1},x_nx_{n-2},w_3,\cdots, w_m)$ are greater than those of $w$.
     If $x_jw_3\in W$,  it is easy to see that the  Anick chain $(x_n, x_n x_{n-1}, x_j, w_3,\cdots,w_m)$ has the greater weight and degree.

 \item[(v)]  $w_2= x_n x_j$ with  $j< n-1$. We can construct a new Anick chain $(x_n, x_n x_{n-1}, w_2, \cdots, w_m)$ whose   weight and degree are greater than those of $w$.
 \end{itemize}
 The above discussion tells us that $w_2=x_nx_{n-1}.$

 Now let us consider $w_3$.  If  $x_{n-1}w_3 \in  W$,    the Anick chain $(x_n, x_nx_{n-1}, x_{n-1}, w_3,\cdots, w_m)$  has weight  and degree greater than those of  $w$.
So it must be $x_n x_{n-1}w_3\in W. $ It follows that  either $w_3=x_j,$ $j<n$ or $w_3=x_n x_{k},$ $k<n-1$.
If  $w_3=x_j$ with $j<n,$  we can construct a new Anick chain $(x_n, x_n x_{n-1}, x_n x_{n-2}, w_3, \cdots, w_m)$ with greater weight  and degree.
  If $w_3= x_n x_{k}$ with $k< n-2$, the Anick chain $(x_n, x_n x_{n-1}, x_n x_{n-2}, w_3, \cdots, w_m)$ has weight  and degree greater than those of $w$. Hence we have $w_3=x_nx_{n-2}$.

Similarly, one can show that
$$w_4= x_n x_{n-3},\cdots, w_{n}= x_n x_1.$$

 Notice that $x_1w_{n+1}\notin W$ for any $w_{n+1}\in X^{*},$ so it must be $x_nx_1w_{n+1}\in W$ which shows that $w_{n+1} \in \{ x_1 , \cdots, x_{n-1} \}.$ If $w_{n+1} \neq x_{n-1}$, the  Anick chain $$(x_n, x_n x_{n-1}, \cdots, x_n x_1, x_{n-1},x_{n-1} x_{n-2}, w_{n+1}, w_{n+2},\cdots, w_m)$$ has the greater  weight  and degree. Hence  $w_{n+1}= x_{n-1}$.

 As the first letter is maximal among all letters in any monomial of $W$, each $w_i$ does not contain $x_n$ for $i\geq n+1.$
 By the induction hypothesis, the segment $(x_{n-1}, w_{n+2},\cdots, w_m)$ of $w$ is of weight $\frac{n(n-1)}{2}$ and its degree is greater than that of any other Anick chain of the form $(x_{n-1},w_2',\cdots,w_p').$
 It follows that $w$ is of weight $\frac{n(n-1)}{2}+n=\frac{n(n+1)}{2}$ and has the maximal degree. Hence (a)(b) holds.


Since $W$ is finite, the Ufnarovski\u\i~graph  $Q_W$    is a finite quiver and the nonexistence of oriented cycles  can be deduced from  (b).  Thus (c) holds.

The statement (d) follows from (c).
\end{Proof}

\begin{Rem}\label{Rem: degree of generator}
 It should be noted that in the proof of Lemma~\ref{Lem: chain structure of Chinese alg}, we  constructed the unique Anick chain of maximal weight and degree $$w_{\mathrm{max}}=(x_n, x_n x_{n-1}, \cdots, x_n x_1, x_{n-1}, x_{n-1} x_{n-2}, \cdots,   x_2, x_2x_1,  x_1)$$
whose degree is equal to $n^2.$ \end{Rem}

\smallskip

Recall that an algebra    is  homologically smooth, if it admits a finite length resolution by finitely generated projective bimodules \cite[Definition 3.1.3]{Ginzburg}.

\begin{Thm}\label{Thm: Chinese algebras}
The Chinese algebra  $A$ of   rank $n\geq 1$ is  homologically smooth and of  global dimension   $\frac{n(n+1)}{2}$.
\end{Thm}

\begin{Proof}
 By Lemma~\ref{Lem: chain structure of Chinese alg},   each $P_i= A\otimes_k kW^{(i-1)} \otimes_k A$ is  a  finitely generated free bimodule and  the two-sided Anick resolution $(P_{*}, d_{*} )$ is of length $\frac{n(n+1)}{2}$. This implies that  $A$ is homologically smooth.

 \medskip

For any left $A$ module $M$,  applying the functor $- \otimes_A M$ to the two-sided Anick resolution yields a free resolution $(P_*\otimes_A M, d_{*}\otimes 1_M) $ of $A\otimes_A M \cong M$ with length $\frac{n(n+1)}{2}$. So we have $\gldim(A)\leq \frac{n(n+1)}{2}$.

On the other hand,  in Lemma~\ref{Lem: chain structure of Chinese alg}, we have shown  that the degree of the generator $w_{max}$ of $P_{\frac{n(n+1)}{2}}$ is greater than any other  generator of the resolution $(P_{*},d_{*}).$ The homogeneity for the relation of $A$ implies that $w_{max}$ cannot be reduced to any other Anick chain. In terms of the proof of Theorem~\ref{Thm: minimal criterion}, the map $1\otimes_{A} d_{\frac{n(n+1)}{2}-1} \otimes_A 1 $ induced by applying functor $k\otimes_A - \otimes_A k$ to the resolution  $(P_{*},d_{*})$ is a zero map. As $\Tor_{m}^{A}(k,k)\cong \mathrm{H}_{m}(k\otimes_A P'_{*} \otimes_A k,1\otimes_A d'_{*} \otimes_A 1)$ for any two-sided projective resolution $(P'_{*}, d'_{*})$ of $A.$ So,  by  the fact that $d_{\frac{n(n+1)}{2}}=0$ and $1\otimes_A d_{\frac{n(n+1)}{2}-1} \otimes_A 1=0,$ we have the isomorphisms  $$\Tor_{\frac{n(n+1)}{2}}^{A}(k,k)\cong k\otimes_A P_{\frac{n(n+1)}{2}} \otimes_A k\cong   k\otimes_A A \otimes kW^{(\frac{n(n+1)}{2}-1)} \otimes A \otimes_A k \cong k W^{(\frac{n(n+1)}{2}-1)} \neq 0.$$
It implies that  the projective dimension of $k$ is as least $\frac{n(n+1)}{2}$. Hence we have $\gldim(A) \geq \frac{n(n+1)}{2}.$  This finishes the proof.

\end{Proof}

\begin{Rem}  The above proof indicates that the two-sided Anick resolution of the Chinese algebra of rank $n\geq 3$ is of minimal length, however, it is  NOT minimal. In fact,  $1_k\ot_A d_3\ot_A 1_k\neq 0$.  \end{Rem}

\begin{Rem}
Although the Chinese algebra $A$ of rank $n\geq 3$ is cubic, it is NOT a $3$-Koszul algebra \cite{Berger, YZ, GMMZ}.
In fact, assume that  $A$ is a $3$-Koszul algebra  generated in degree $1$, then, by definition, $\Tor_{i}^{A}(k,k)$ has a basis in degree $\delta(i),$ where
$$\delta(i)=\left\{ \begin{array}{lc} \frac{3i}{2} & \text{if $i$ is even} \\
\frac{3(i-1)}{2}+1 & \text{if $i$ is odd} \end{array}\right.  .$$
We have proved that $\Tor_{\frac{n(n+1)}{2}}^{A}(k,k)=k\{(x_n,x_nx_{n-1},\cdots, x_nx_1,x_{n-1}, \cdots, x_1)\}$ in Theorem~\ref{Thm: Chinese algebras} and the degree of the unique generator of $\Tor_{\frac{n(n+1)}{2}}^{A}(k,k)$ is $n^2$ by Remark~\ref{Rem: degree of generator}  which is not equal to $ \delta(\frac{n(n+1)}{2})$,  thus giving  a contradiction.

\end{Rem}

 In a forthcoming paper, we will try to construct a minimal two-sided projective resolution of the Chinese algebra and compute its Hochschild cohomology.

\bigskip

\section{A Koszul algebra whose Anick resolution is not minimal}

  In order to answer two questions from the book \cite{PP}, N.~Iyudu and S.~Shkarin \cite{IS} classified Hilbert series of Koszul algebras with three generators and three relations. They introduce a new  Koszul algebra $$A=k \langle x,y,z \ | \ x^2+yx ,xz,zy\rangle.$$
V.~Dotsenko and S.~R.~Chowdhury \cite{Dots} calculated the bar homology $\Tor_{*}^{A}(k,k)$ of $A$    through the (one-sided) Anick resolution. In this section, we will use algebraic Morse theory to construct the two-sided minimal resolution of $A$ from the two-sided Anick resolution which itself is not minimal.

\begin{Lem}\cite{Dots}
Let $A$ be as above. The Gr\"{o}bner-Shirshov basis  for the ideal $I$ with respect to the left length-lexicographic order generated by $x>y>z$ is $$\mathcal{G}=\{ xy^{k}x + y^{k+1}x ,xz,zy\}_{k\geq 0}.$$
\end{Lem}

 It follows that  the reduced monomial generators of $\Tip(I)$ is
$$W=\Tip(\mathcal{G})=\{ xy^{k}x ,xz,zy\}_{k\geq 0}. $$
It produces the list of all Anick chains as follows.
\begin{itemize}
\item[(a)] the set of $(-1)$-chain $W^{(-1)}$ is a singleton;
\item[(b)] the set of $0$-chains $W^{(0)}$ consists of $x,y,z;$
\item[(c)] the set of $1$-chains $W^{(1)}$ consists of  $(x,y^{k}x),(x,z),(z,y)$ with $k\geq 0$;
\item[(d)] for $n\geq 2,$ the set of $n$-chains $W^{(n)}$ consists of $$(x,y^{k_1}x,\cdots, y^{k_n}x ), (x,y^{k_1}x,\cdots,y^{k_{n-1}}x,z), (x,y^{k_1}x,\cdots,y^{k_{n-2}}x,z,y) $$
    with $k_1,\cdots , k_{n} \geq 0.$
\end{itemize}
By Theorems~\ref{Thm: main result of algebraic Morse theory} and \ref{thm:skoeldberg}, we obtain a free bimodule resolution $(P_{*}, d_{*} )$ of $A$ with
$$P_{i}=A \otimes k W^{(i-1)} \otimes A, i\geq 0.$$

\begin{Prop}\label{Prop: differential of the quadratic algebra}
The differential of the two-sided resolution is given by the following.
\begin{itemize}
 \item[(a)] $d_{0}(1\otimes 1)=1;$
 \item[(b)] $d_{1}(1\otimes a \otimes 1)=a\otimes 1- 1\otimes a, $  $a=x,y,z;$
 \item[(c)] $ d_{2} (1\otimes (x,y^{k}x) \otimes 1)= (xy^{k} + y^{k+1}) \otimes x \otimes 1 + 1\otimes x \otimes y^{k}x +1\otimes y \otimes y^{k}x + \sum\limits_{i=1}^{k} (xy^{i-1} +y^{i} )\otimes y \otimes y^{k-i}x , \\
     d_{2} (1\otimes (x,z) \otimes 1) =x\otimes z \otimes 1 +1\otimes x \otimes z, \\ d_{2}(1\otimes (z,y) \otimes 1) = z\otimes y \otimes 1 + 1\otimes z \otimes y;$
 \item[(d)] for $n\geq 3,$

  $d_{n+1}(1\otimes (x,y^{k_1}x,\cdots, y^{k_n}x)\otimes 1) $
\begin{equation}\label{Eq: d_n+1 of Anick resolution 1}\begin{split}
 = & (xy^{k_1} + y^{k_1+1} )\otimes (x,y^{k_2}x,\cdots, y^{k_n}x)\otimes 1 +  (-1)^{n+1}\otimes (x,y^{k_1}x,\cdots, y^{k_{n-1}}x)\otimes y^{k_n}x \\
& +  \sum_{i=1}^{n-1}(-1)^{i}  \otimes (x,y^{k_1}x,\cdots, y^{k_{i-1}}x , y^{k_{i}+k_{i+1}+1}x, y^{k_{i+2}}x,\cdots, y^{k_n}x) \otimes 1 ,
 \end{split}\end{equation}

$ d_{n+1}(1\otimes (x,y^{k_1}x,\cdots, y^{k_{n-1}}x, z)\otimes 1)$
 \begin{equation} \label{Eq: d_n+1 of Anick resolution 2}\begin{split}
 = &(xy^{k_1} +y^{k_1+1}) \otimes (x,y^{k_2}x,\cdots, y^{k_{n-1}}x,z)\otimes 1 +  (-1)^{n+1}\otimes (x,y^{k_1}x,\cdots, y^{k_{n-1}}x)\otimes z \\
 & + \sum_{i=1}^{n-2}(-1)^{i} \otimes (x,y^{k_1}x,\cdots, y^{k_{i-1}}x , y^{k_{i}+k_{i+1}+1}x, y^{k_{i+2}}x,\cdots, y^{k_{n-1}}x,z) \otimes 1 ,
 \end{split} \end{equation}

 $ d_{n+1}(1\otimes (x,y^{k_1}x,\cdots, y^{k_{n-2}}x, z, y)\otimes 1)$
 \begin{equation}\label{Eq: d_n+1 of Anick resolution 3}\begin{split}
 = &(xy^{k_1} +y^{k_1+1}) \otimes (x,y^{k_2}x,\cdots, y^{k_{n-2}}x,z,y)\otimes 1 +  (-1)^{n+1}\otimes (x,y^{k_1}x,\cdots, y^{k_{n-2}}x,z)\otimes y \\
 & + \sum_{i=1}^{n-3}(-1)^{i} \otimes (x,y^{k_1}x,\cdots, y^{k_{i-1}}x , y^{k_{i}+k_{i+1}+1}x, y^{k_{i+2}}x,\cdots, y^{k_{n-2}}x,z,y) \otimes 1 .
 \end{split} \end{equation}

\end{itemize}
\end{Prop}

\begin{Proof}
At the risk of being repetitive, again, as we are working with Theorem~\ref{Thm: main result of algebraic Morse theory}, the differential of $(P_{*},d_{*})$ is determined by all the zigzag paths between two critical vertices in the quiver $\overline{Q_B}^\M$ (cf. Section 3  and Remark \ref{Def:new-weighted-quiver-overline}).

The maps $d_{0}$ and $d_{1}$ are easy. Now let us consider $d_{2}(1\otimes (x,y^{k}x) \otimes 1),$ or equivalently, all the zigzag paths in $\overline{Q_B}^\M$ with starting vertex being $(x,y^{k}x)$ and target lying in $W^{(0)}.$ Let $p=\alpha_1\cdots \alpha_m$ be such a zigzag path. Notice that $\alpha_1$ and $\alpha_m$ must be the thick arrow. We discuss $\alpha_1$ in three cases (the notations $d_n^{j}$ below are given in Section 4).
\begin{itemize}
\item[(i)] $\alpha_1=d_2^2:(x,y^{k}x)  \overset{1\otimes y^{k}x}{\to} x .  $ As $x\in W^{(0)}$ is a critical vertex, there is no dotted arrow with $x$ as its starting vertex. Thus the zigzag path $p=\alpha_1$ and it gives the term $ 1\otimes x \otimes y^{k}x$ of $d_{2}(1\otimes (x,y^{k}x) \otimes 1  );$
\item[(ii)] $\alpha_1=d_2^0: (x,y^{k}x)  \overset{x\otimes 1}{\to}  y^{k}x .$ Then $\alpha_2=d_2^{-1}: y^{k}x  \overset{1\otimes 1}{\dashrightarrow}  (y,y^{k-1}x) $ is the unique dotted arrow with $y^{k}x  $ as its starting vertex. It follows that $\alpha_3$ is either $d_2^0$ or $d_2^2$ as the arrow $d_2^{1} : (y,y^{k-1}x)  \overset{-1\otimes 1}{\to}   y^{k}x  $  does not exist in $\overline{Q_B}^\M.$
    \begin{itemize}
    \item[(ii-1)] if $\alpha_3=d_2^2:  (y,y^{k-1}x)  \overset{1\otimes y^{k-1}x}{\to} y ,$  then the zigzag path is $p=\alpha_1\alpha_2\alpha_3$ and it gives the term $ x \otimes y \otimes y^{k-1}x$ of $d_{2}(1\otimes (x,y^{k}x) \otimes 1  );$
    \item[(ii-2)] if $\alpha_3=d_2^0 : (y,y^{k-1}x)  \overset{y\otimes 1}{\to}   y^{k-1}x , $ then the discussion of the zigzag path $p'=\alpha_3\cdots \alpha_m$ is  similar as that in the case (ii).
    \end{itemize}
    It follows   by induction  that all the zigzag paths with $\alpha_1=d_2^0$ as their first thick arrows give the terms
    $$\sum_{i=1}^{k} xy^{i-1} \otimes y \otimes y^{k-i}x + xy^{k} \otimes x\otimes 1$$
     of $d_{2}(1\otimes (x,y^{k}x) \otimes 1  );$
    \item[(iii)] $\alpha_1=d_{2}^{1}: (x,y^{k}x)  \overset{1\otimes 1}{\to} y^{k+1}x.$ This case is also similar to the case (ii) and one can show that the zigzag paths with first thick arrows being  $\alpha_1=d_2^1$  give the terms $$\sum_{i=1}^{k+1} y^{i-1} \otimes y \otimes y^{k-i+1}x + y^{k+1}\otimes x \otimes 1 $$
        of $d_{2}(1\otimes (x,y^{k}x) \otimes 1).$
\end{itemize}
Combining the results above, we obtain the  expression of $d_2$ on $1\otimes (x,y^{k}x) \otimes 1.$ The formulas of
$d_{2}(1\otimes (x,z) \otimes 1)$ and $d_{2}(1\otimes (z,y) \otimes 1) $ can be proved similarly.

For $n\geq 2,$ let $p=\alpha_1\cdots \alpha_m$ be a zigzag path with starting vertex being $(x,y^{k_1}x,\cdots, y^{k_n}x)$ and target lying in $W^{(n-1)}.$ We discuss $\alpha_1$ in four cases.
\begin{itemize}
\item[(i)] $\alpha_1=d_{n+1}^{n+1}: (x,y^{k_1}x,\cdots, y^{k_n}x) \overset{(-1)^{n+1}\otimes y^{k_n}x}{\to}  (x,y^{k_1}x,\cdots, y^{k_{n-1}}x)\in W^{(n-1)}.$ We have $p=\alpha_1$ and this zigzag path gives the term $(-1)^{n+1}\otimes (x,y^{k_1}x,\cdots, y^{k_{n-1}}x) \otimes y^{k_n}x$ of $d_{n+1}(1\otimes  (x,y^{k_1}x,\cdots, y^{k_n}x) \otimes 1);$
\item[(ii)] $\alpha_1=d_{n+1}^{i+1}:(x,y^{k_1}x,\cdots, y^{k_n}x) \overset{(-1)^{i}\otimes 1}{\to}  (x,y^{k_1}x,\cdots, y^{k_{i-1}}x , y^{k_{i}+k_{i+1}+1}x,y^{k_{i+2}}x,\cdots, y^{k_n}x)\in W^{(n-1)}$ with $1\leq i \leq n-1.$ It follows that  $p=\alpha_1$ and these zigzag paths give the terms $(-1)^{i}\otimes (x,y^{k_1}x,\cdots, y^{k_{i-1}}x , y^{k_{i}+k_{i+1}+1}x,y^{k_{i+2}}x,\cdots, y^{k_n}x) \otimes 1$ of $d_{n+1}(1\otimes  (x,y^{k_1}x,\cdots, y^{k_n}x) \otimes 1);$
\item[(iii)] $\alpha_1=d_{n+1}^{0}: (x,y^{k_1}x,\cdots, y^{k_n}x) \overset{x\otimes 1}{\to} (y^{k_1}x,\cdots, y^{k_n}x) .$ Then $\alpha_2=d_{n+1}^{-1} : (y^{k_1}x,\cdots, y^{k_n}x) \overset{1\otimes 1}{\dashrightarrow}(y,y^{k_1-1}x,\cdots, y^{k_n}x)  .$ If $\alpha_3$ is of the form $d_{n+1}^{i}:(y,y^{k_1-1}x,\cdots, y^{k_n}x) \to (y, y^{k_{1}-1}x,\cdots) $ with $2\leq i \leq n+1,$ there is no dotted arrow with starting vertex being $(y, y^{k_{1}-1}x,\cdots). $ Meanwhile, $\alpha_3$ could not be $d_{n+1}^{1}.$ Thus $\alpha_3$ can only be $d_{n+1}^{0}:(y,y^{k_1-1}x,\cdots, y^{k_n}x) \overset{y\otimes 1}{\to} (y^{k_1-1}x,\cdots, y^{k_n}x).$ Repeating the discussion in this case for
    $p=\alpha_3\cdots \alpha_m$ and using the induction on it, we can see that the zigzag path in this case is $p=d_{n+1}^{0} d_{n+1}^{-1} d_{n+1}^{0}\cdots d_{n+1}^{-1} d_{n+1}^{0}$ and it gives the term $ xy^{k_1} \otimes (x,y^{k_2}x,\cdots, y^{k_n}x) \otimes 1 $ of $d_{n+1}(1\otimes  (x,y^{k_1}x,\cdots, y^{k_n}x) \otimes 1);$
\item[(iv)] $\alpha_1=d_{n+1}^{1}: (x,y^{k_1}x,\cdots, y^{k_n}x) \overset{1\otimes 1}{\to} (y^{k_1+1}x,\cdots, y^{k_n}x) .$ Similar to the discussion in case (iii), one can show that $p=d_{n+1}^{1} d_{n+1}^{-1}d_{n+1}^{0}d_{n+1}^{-1}\cdots d_{n+1}^{-1}d_{n+1}^{0}$ and it gives the term $ y^{k_1+1} \otimes (x,y^{k_2}x,\cdots, y^{k_n}x) \otimes 1$ of $d_{n+1}(1\otimes (x,y^{k_1}x,\cdots, y^{k_n}x) \otimes 1).$
\end{itemize}
The above discussions prove the formula of $d_{n+1}(1\otimes  (x,y^{k_1}x,\cdots, y^{k_n}x) \otimes 1).$ The remaining formulas can be proved  similarly.
\end{Proof}

Denote the new weighted quiver  with respect to the Anick resolution $(P_{*},d_{*})$ by $Q_{P}$ (cf. Section 3). Its vertex set is the set of all the Anick chains and arrow set is determined by the differential in Proposition~\ref{Prop: differential of the quadratic algebra}. We define a full subquiver $\M$ of $Q_{P}$ as the union of the following sets of arrows.

$$\begin{array}{c}
\{  (x,x,y^{\ell_1}x,\cdots , y^{\ell_n} x) \overset{-1\otimes 1}{\to} (x,y^{\ell_1+1}x,\cdots , y^{\ell_n} x) ; n\geq 1, \ell_1,\cdots, \ell_n\geq 0 \} \\
 \{  (x,x,y^{\ell_1}x,\cdots , y^{\ell_n} x,z) \overset{-1\otimes 1}{\to} (x,y^{\ell_1+1}x,\cdots , y^{\ell_n} x,z) ; n\geq 1, \ell_1,\cdots, \ell_n\geq 0 \}\\
  \{  (x,x,y^{\ell_1}x,\cdots , y^{\ell_n} x,z,y) \overset{-1\otimes 1}{\to} (x,y^{\ell_1+1}x,\cdots , y^{\ell_n} x,z,y) ; n\geq 1, \ell_1,\cdots, \ell_n\geq 0 \}
\end{array}$$
It is easy to see that each arrow in $\M$ has weight $-1\otimes 1$ and two different arrows in $\M$ do not share the endpoints. Thus we have the following.
\begin{Lem}
$\M$ is a partial matching.
\end{Lem}

In the quiver $Q_{P}^{\M},$ two endpoints of a dotted arrow share the degree (number of variables) and for a thick arrow, the degree of its starting vertex  is greater or equal to that of its target. Let $p=\alpha_1\cdots\alpha_m$ be a  zigzag path in $Q_{P}^{\M}$. Without loss of generality, we may assume that all arrows in $p$ preserve the degree of vertices. Consider the segment $v_i \overset{\alpha_i}{\dashrightarrow} v_{i+1} \overset{\alpha_{i+1}}{\to} v_{i+2}$ of $p.$
 By the construction of $\M,$  $\alpha_i$ has   the form $\alpha_i: (x,y^{\ell_1+1}x, \cdots ) \dashrightarrow (x,x,y^{\ell_1}x, \cdots ).$
  Note that we are concerned with the quiver  $Q_{P}^{\M}$ rather than $Q_{B}^{\M},$ the thick arrow $\alpha_{i+1}$ is given by the differential of the two-sided Anick resolution  \eqref{Eq: d_n+1 of Anick resolution 1} \eqref{Eq: d_n+1 of Anick resolution 2} \eqref{Eq: d_n+1 of Anick resolution 3}. Assume that $\alpha_{i+1}$ corresponds to the term $(x+y) \otimes (x,y^{\ell_1}x,\cdots ) \otimes 1$ of $d(1\otimes (x,x,y^{\ell_1}x, \cdots ) \otimes 1 ),$ that is, the terminal vertex $v_{i+2}$ of $\alpha_{i+1} $ is $(x,y^{\ell_1}x,\cdots ) ,$ it follows that two endpoints of $\alpha_{i+1}$ do not share the degree. According to the construction of $Q_{P}^{\M}$, we have $v_{i+2}\neq v_{i}.$ Therefore, under our assumption that $\alpha_{i+1}$ preserves the degree of vertices, $v_{i+2} $ has the form $(x,x,y^{\ell_1}x, \cdots).$
By the choice of the arrows of $\M$ there is no dotted arrow with $v_{i+2}$ as its starting vertex. Thus the zigzag path $p$  is of finite length. According to Proposition~\ref{Prop: sufficient condition for Morse matching}, we have
\begin{Prop}
$\M$ is a Morse matching.
\end{Prop}

Obviously, the set of critical vertices $\V^{\M}$ consists of the Anick chains without being $(x,x,y^{\ell_1}x,\cdots ) $ and $(x,y^{\ell_1+1}x,\cdots )$ with $\ell_1\geq 0.$ This observation filters out most vertices and we list the vertices remaining as follows.
$$\V^{\M}_4 = \{ (x,x,z,y)\} ,\V_{3}^{\M}=\{ (x,x,z), (x,z,y)\} , \V_{2}^{\M}=\{ (x,x), (x,z),(z,y)\}, \V_{1}^{\M}=\{x,y,z\}, \V_{0}^{\M}=\{*\} .  $$
By Theorem~\ref{Thm: main result of algebraic Morse theory}, we obtain the two-sided resolution of $A$
$$0 \to A\otimes k\V^{\M}_4  \otimes A \overset{d_4'}{\to} A\otimes k\V^{\M}_3 \otimes A \overset{d_3'}{\to} A\otimes k \V^{\M}_2  \otimes A \overset{d_2'}{\to} A \otimes k \V^{\M}_1  \otimes A \overset{d_1'}{\to} A \otimes A \overset{\epsilon}{\to} A,$$
The differential listed below can be easily calculated by the enumeration of zigzag paths of $Q_{P}^{\M}.$
$$\begin{array}{l}
d_4'(1\otimes (x,x,z,y) \otimes 1) = (x+y)\otimes (x,z,y) \otimes 1 + 1\otimes(x,x,z) \otimes y ,\\
d_3'(1\otimes(x,x,z)\otimes 1)= (x+y)\otimes (x,z)\otimes 1 - 1\otimes(x,x) \otimes z, \\
d_{3}'(1\otimes (x,z,y)\otimes 1)= x\otimes(z,y) \otimes 1 - 1\otimes (x,z) \otimes y,\\
d_{2}'(1\otimes (x,x) \otimes 1)=(x+y)\otimes x \otimes 1 +1\otimes x \otimes x + 1\otimes y \otimes x,\\
d_{2}'(1\otimes (x,z) \otimes 1)=x\otimes z \otimes 1 +1\otimes x \otimes z,\\
d_{2}'(1\otimes (z,y) \otimes 1)=z\otimes y \otimes 1 +1\otimes z \otimes y,\\
d_1'(1\otimes a \otimes 1)= a\otimes 1- 1\otimes a, a=x,y,z,
\end{array}$$
$\epsilon$ is the multiplication. Hence this resolution is minimal.  Note that by tensoring $k$ over $A$ from the right we get a minimal one-sided resolution which is smaller than the one-sided Anick resolution obtained in \cite{Dots}.

 In a forthcoming paper, the Hochschild cohomology groups of $A$ as well as its Gerstenhaber algebra structure will be determined explicitly.

\bigskip

\textbf{Acknowledgements:}   The    authors were   supported by the National Natural Science Foundation of
China (No. 12071137, No. 12031014), by Key Laboratory of MEA (Ministry of Education), by Shanghai Key
Laboratory of PMMP (No. 22DZ2229014), and by Fundamental Research Funds for the Central
Universities.

It is our great pleasure to thank Severin Barmeier who has read this paper carefully and  has  given us many suggestions to improve it.

After the submission of this paper, we found a paper by H.~Alhussein \cite{Alh22} who also used algebraic Morse theory to compute Anick resolutions and Hochschild cohomologies of the Chinese monoid algebra for small ranks.

\end{document}